\documentclass[12pt,a4paper]{amsart}

\usepackage{hyperref} 
\usepackage{amssymb}
\usepackage{amsrefs}
\usepackage[all]{xy}

\allowdisplaybreaks

\title[On derived equivalence classification of gentle algebras]%
  {On derived equivalence classification of gentle two-cycle algebras}

\author{Grzegorz Bobi\'nski}

\author{Piotr Malicki (Toru\'n)}

\address{Faculty of Mathematics and Computer Science \\ Nicolaus
Copernicus University \\ ul.~Chopina 12/18 \\ 87-100 Toru\'n \\
Poland}

\email{gregbob@mat.uni.torun.pl}

\email{pmalicki@mat.uni.torun.pl}

\keywords{derived category, gentle algebra, tilting-cotilting
equi\-va\-len\-ce}

\subjclass[2000]{18E30, 16G20}

\newtheorem{theo}{Theorem}
\newtheorem*{conj}{Conjecture}

\newtheorem{claim}{Claim}[section]

\newtheorem{coro}[claim]{Corollary}
\newtheorem{lemm}[claim]{Lemma}
\newtheorem{prop}[claim]{Proposition}

\newcounter{move}[section]
\newtheorem{step}[move]{Step}

\theoremstyle{remark}

\DeclareMathOperator{\id}{id}%
\DeclareMathOperator{\op}{op}%
\DeclareMathOperator{\pd}{pd}%
\DeclareMathOperator{\rk}{rk}%
\DeclareMathOperator{\End}{End}%
\DeclareMathOperator{\Ext}{Ext}%

\newcommand{\bbA}{\mathbb{A}}
\newcommand{\bbN}{\mathbb{N}}

\newcommand{\bbZ}{\mathbb{Z}}

\newcommand{\calA}{\mathcal{A}}
\newcommand{\calB}{\mathcal{B}}
\newcommand{\calD}{\mathcal{D}}

\newcommand{\vertexD}[1]{\bullet \save*+!D{\scriptstyle #1} \restore}
\newcommand{\vertexL}[1]{\bullet \save*+!L{\scriptstyle #1} \restore}
\newcommand{\vertexR}[1]{\bullet \save*+!R{\scriptstyle #1} \restore}
\newcommand{\vertexU}[1]{\bullet \save*+!U{\scriptstyle #1} \restore}

\newlength{\branch}
\begin{document}

\begin{abstract}
We classify, up to derived (equivalently, tilting-co\-til\-ting)
equivalence all nondegenerate gentle two-cycle algebras. We also
give a partial classification and formulate a conjecture in the
degenerate case.
\end{abstract}

\maketitle

\section*{Introduction and the main result}

Throughout the paper $k$ denotes a fixed algebraically closed
field. By an algebra we mean a finite dimensional basic connected
$k$-algebra and by a module a finite dimensional left module. By
$\bbZ$, $\bbN$, and $\bbN_+$, we denote the sets of integers,
nonnegative integers, and positive integers, respectively.
Finally, if $i, j \in \bbZ$, then $[i, j] = \{ l \in \bbZ \mid i
\leq l \leq j \}$.

With an algebra $\Lambda$ we may associate its bounded derived
category $\calD^b (\Lambda)$ (in the sense of
Verdier~\cite{Ver1977}) of bounded complexes of $\Lambda$-modules,
which has a structure of a triangulated category
(see~\cite{Hap1988}). The bounded derived category is an important
homological invariant of the module category of an algebra and
attracts a lot of interest (see for example~\cites{AssSko1988,
BarLen2003, Gei2002, Hap1987, Hap1991, Kel1994, Ric1989a,
Ric1989b}). In particular, the derived equivalence classes of
algebras have been investigated (see for example~\cites{Asa1999,
BiaHolSko2003, BocHolSko2006, Bru2001, Hol1999}), where two
algebras are said to be derived equivalent if their bounded
derived categories are equivalent as triangulated categories.

A handy way of proving a derived equivalence between algebras
$\Lambda$ and $\Lambda'$ is a construction of a (co)tilting
$\Lambda$-module $T$ such that $\Lambda'$ is (isomorphic to) the
opposite of the endomorphism algebra of $T$. Here a
$\Lambda$-module $T$ is called (co)tilting if $\pd_\Lambda T \leq
1$ ($\id_\Lambda T \leq 1$, respectively), $\Ext_\Lambda^1 (T, T)
= 0$, and $T$ is a direct sum of precisely $\rk K_0 (\Lambda)$
pairwise nonisomorphic indecomposable $\Lambda$-modules, where
$K_0 (\Lambda)$ denotes the Grothendieck group of the category of
$\Lambda$-modules. The transitive closure of the relation defined
in this way is called tilting-cotilting equivalence. For many
classes of algebras tilting-cotilting equivalence and derived
equivalence coincide.

Results of this type have been obtained for gentle algebras,
introduced by Assem and Skowro\'nski in~\cite{AssSko1987} (see
Section~\ref{sect_prem} for a precise definition), which form an
important subclass of the class of special biserial algebras in
the sense of~\cite{SkoWas1983}. We note that a
representation-infinite algebra admits a simply connected Galois
covering whose every finite convex subcategory is
representation-finite if and only if it is a special biserial
algebra and its simply connected Galois covering is the repetitive
category of union of a countable chain of gentle tree algebras
(see~\cite{PogSko1991}). The class of algebras derived equivalent
to a hereditary algebra of Dynkin type $\bbA_n$ for some $n \in
\bbN_+$ coincides with the class of algebras tilting-cotilting
equivalent to a hereditary algebra of type $\bbA_n$ and consists
of the gentle algebras whose Gabriel quivers have $n$ vertices and
$n - 1$ arrows (see~\cite{AssHap1981}). Moreover, for a given $n$
all such algebras form one derived equivalence class. Similarly,
the class of algebras derived equivalent to a hereditary algebra
of Euclidean type $\tilde{\bbA}_n$ for some $n \in \bbN_+$
coincides with the class of algebras tilting-cotilting equivalent
to a hereditary algebra of type $\tilde{\bbA}_n$ and consists of
the gentle algebras whose Gabriel quivers have $n$ vertices and
$n$ arrows and which satisfy the so-called clock condition on the
unique cycle. In this case, there are exactly $\lfloor \frac{n}{2}
\rfloor$ derived (equivalently, tilting-cotilting) equivalence
classes for a given $n$. The algebras with the same numbers of
vertices and arrows in the Gabriel quiver equal are called
one-cycle algebras. The remaining gentle one-cycle algebras form
the class of derived discrete algebras which are not derived
(equivalently, tilting-cotilting) equivalent to a hereditary
algebra of Dynkin type (see~\cite{Vos2001}). The derived
equivalence classes of these algebras were described
in~\cite{BobGeiSko2004}.

The aim of this paper it to extend the above classification to the
class of gentle two-cycle algebras, where we call an algebra a
two-cycle algebra if the number of arrows in the Gabriel quiver
exceeds the number of vertices by one. An additional motivation
for this research is the fact proved by Schr\"oer and Zimmermann
in~\cite{SchZim2003} saying that the gentle algebras are closed
under derived equivalences. Moreover, for the gentle algebras the
numbers of vertices and arrows in the Gabriel quiver are derived
invariants (see~\cite{AAGei2006}*{Corollary~15}). However, we
obtain a full classification only for nondegenerate gentle
two-cycle algebras, where we call a gentle two-cycle algebra
$\Lambda$ nondegenerate if $\sum_{(n, m) \in \bbN \times \bbN}
\phi_\Lambda (n, m) = 3$ and $\phi_\Lambda : \bbN \times \bbN \to
\bbN$ is the derived invariant introduced by Avella-Alaminos and
Geiss in \cite{AAGei2006} (see Section~\ref{sect_min}). For the
remaining gentle two-cycle algebras $\Lambda$, which we call
degenerate, we have $\sum_{(n, m) \in \bbN \times \bbN}
\phi_\Lambda (n, m) = 1$. Obviously, the both classes of gentle
two-cycle algebras are closed under derived (hence also
tilting-cotilting) equivalences.

Before formulating the main results of the paper we define the
following families of algebras.
\begin{itemize}

\item
$\Lambda_0 (p, r)$ for $p \in \bbN_+$ and $r \in [0, p - 1]$ is
the algebra of the quiver
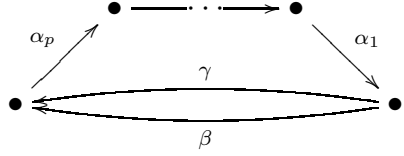
\[
\xymatrix{& \bullet \ar[rr]|{\textstyle \cdots} & & \bullet
\ar[rd]^{\alpha_1} \\ \bullet \ar[ru]^{\alpha_p} & & & & \bullet
\ar@/_/[llll]_\gamma \ar@/^/[llll]^\beta}
\]
bound by $\alpha_p \beta$, $\alpha_i \alpha_{i + 1}$ for $i \in
[1, r]$ and $\gamma \alpha_1$.

\item
$\Lambda_0' (p, r)$ for $p \in \bbN_+$ and $r \in [0, p - 1]$ is
the algebra of the quiver
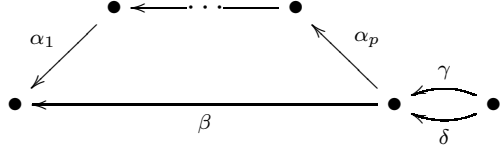
\[
\xymatrix{& \bullet \ar[ld]_{\alpha_1} & & \bullet
\ar[ll]|{\textstyle \cdots} \\ \bullet & & & & \bullet
\ar[llll]^\beta \ar[lu]_{\alpha_p} & \bullet \ar@/_/[l]_\gamma
\ar@/^/[l]^\delta}
\]
bound by $\alpha_i \alpha_{i + 1}$ for $i \in [1, r]$, $\alpha_p
\gamma$, and $\beta \delta$.

\item
$\Lambda_1 (p_1, p_2, p_3, p_4, r_1)$ for $p_1, p_2 \in \bbN_+$,
$p_3, p_4 \in \bbN$, and $r_1 \in [0, p_1 - 1]$, such that $p_2 +
p_3 \geq 2$ and $p_4 + r_1 \geq 1$, is the algebra of the quiver
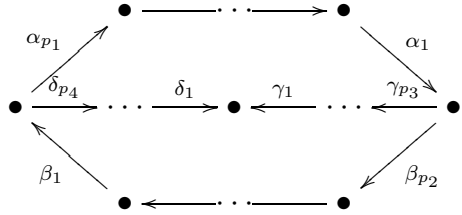
\[
\xymatrix{& \bullet \ar[rr]|{\textstyle \cdots} & & \bullet
\ar[rd]^{\alpha_1} \\ \bullet \ar[ru]^{\alpha_{p_1}}
\ar[r]^-{\delta_{p_4}} & \cdots \ar[r]^-{\delta_1} & \bullet &
\cdots \ar[l]_-{\gamma_1} & \bullet \ar[l]_-{\gamma_{p_3}}
\ar[ld]^{\beta_{p_2}} \\ & \bullet \ar[lu]^{\beta_1} & & \bullet
\ar[ll]|{\textstyle \cdots}}
\]
bound by $\alpha_i \alpha_{i + 1}$ for $i \in [p_1 - r_1, p_1 - 1
]$, $\alpha_{p_1} \beta_1$, $\beta_i \beta_{i + 1}$ for $i \in [1,
p_2 - 1]$, and $\beta_{p_2} \alpha_1$.

\item
$\Lambda_2 (p_1, p_2, p_3, r_1, r_2)$ for $p_1, p_2 \in \bbN_+$,
$p_3 \in \bbN$, $r_1 \in [0, p_1 - 1]$, and $r_2 \in [0, p_2 -
1]$, such that $p_3 + r_1 + r_2 \geq 1$, is the algebra of the
quiver
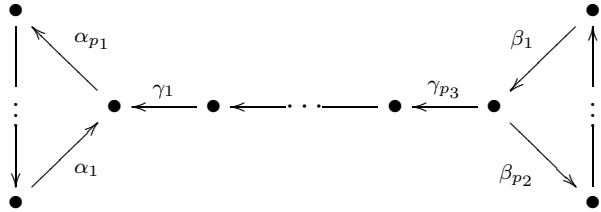
\[
\xymatrix{\bullet \ar[dd]|{\textstyle \vdots} & & & & & & \bullet
\ar[ld]_{\beta_1} \\ & \bullet \ar[lu]_{\alpha_{p_1}} & \bullet
\ar[l]_{\gamma_1} & & \bullet \ar[ll]|{\textstyle \cdots} &
\bullet \ar[l]_{\gamma_{p_3}} \ar[rd]_{\beta_{p_2}} \\ \bullet
\ar[ru]_{\alpha_1} & & & & & & \bullet \ar[uu]|{\textstyle
\vdots}}
\]
bound by $\alpha_i \alpha_{i + 1}$ for $i \in [p_1 - r_1, p_1 -
1]$, $\alpha_{p_1} \alpha_1$, $\beta_i \beta_{i + 1}$ for $i \in
[p_2 - r_2, p_2 - 1]$, and $\beta_{p_2} \beta_1$.

\end{itemize}

The main results of the paper are the following.

\begin{theo} \label{theo_nondeg}
If $\Lambda$ is a nondegenerate gentle two-cycle algebra, then
$\Lambda$ is derived \textup{(}equivalently,
tilting-cotilting\textup{)} equivalent to one of the following
algebras:
\begin{itemize}

\item
$\Lambda_1 (p_1, p_2, p_3, p_4, r_1)$ for some $p_1, p_2 \in
\bbN_+$, $p_3, p_4 \in \bbN$, and $r_1 \in [0, p_1 - 1]$, such
that $p_2 + p_3 \geq 2$, $p_4 + r_1 \geq 1$, and either $p_3
> p_4$ or $p_3 = p_4$ and $p_2 > r_1$,

\item
$\Lambda_2 (p_1, p_2, p_3, r_1, r_2)$ for some $p_1, p_2 \in
\bbN_+$, $p_3 \in \bbN$, $r_1 \in [0, p_1 - 1]$, and $r_2 \in [0,
p_2 - 1]$, such that $p_3 + r_1 + r_2 \geq 1$ and either $p_1 >
p_2$ or $p_1 = p_2$ and $r_1 \geq r_2$.

\end{itemize}
Moreover, different algebras from the above list are not derived
\textup{(}equivalently, tilting-cotilting\textup{)} equivalent.
\end{theo}

\begin{theo} \label{theo_deg}
If $\Lambda$ is a degenerate gentle two-cycle algebra, then
$\Lambda$ is derived \textup{(}equivalently,
tilting-cotilting\textup{)} equivalent to one of the following
algebras:
\begin{itemize}

\item
$\Lambda_0 (p, r)$ for some $p \in \bbN_+$ and $r \in [0, p - 1]$,

\item
$\Lambda_0' (p, 0)$ for some $p \in \bbN_+$.
\end{itemize}
\end{theo}

Moreover, we have the following conjecture concerning the
minimality of the list in the above theorem.

\begin{conj}
Different algebras from the list in Theorem~\ref{theo_deg} are not
derived \textup{(}equivalently, tilting-cotilting\textup{)}
equivalent.
\end{conj}

Obviously, if $p_1, p_2 \in \bbN_+$, $r_1 \in [0, p_1 - 1]$, $r_2
\in [0, p_2 - 1]$, and $p_1 \neq p_2$, then $\Lambda_0 (p_1, r_1)$
and $\Lambda_0 (p_2, r_2)$ ($\Lambda_0' (p_1, 0)$ and $\Lambda_0'
(p_2, 0)$, respectively) are not derived equivalent. Similarly, if
$p_1, p_2 \in \bbN_+$, $r_1 \in [0, p_1 - 1]$, and $p_1 \neq p_2 +
1$, then $\Lambda_0 (p_1, r_1)$ and $\Lambda_0' (p_2, 0)$ are not
derived equivalent. Thus it is enough to prove that $\Lambda_0 (p
+ 1, 0)$, \ldots, $\Lambda_0 (p + 1, p)$ and $\Lambda_0' (p, 0)$
are pairwise not derived equivalent for a fixed $p \in \bbN_+$. It
follows easily by investigating the Euler quadratic forms that
$\Lambda_0 (p + 1, r_1)$ and $\Lambda (p + 1, r_2)$ ($\Lambda_0'
(p, 0)$ and $\Lambda_0 (p + 1, r_2)$) are not derived equivalent
if $r_1 \not \equiv r_2 \pmod{2}$ ($r_2 \equiv 0 \pmod{2}$,
respectively).

The paper is organized as follows. In Section~\ref{sect_prem} we
first present basic definitions, then describe main tools used in
order to reduce an arbitrary gentle two-cycle algebra to one of
the algebras listed in Theorems~\ref{theo_nondeg}
and~\ref{theo_deg}: passing to the opposite algebra, (generalized)
APR-(co)reflections, and HW-(co)reflections. Finally, we describe
in Section~\ref{sect_prem} an operation of shifting relations
being a basic application of the above operations, and investigate
two particular families of gentle two-cycle algebras. In
Section~\ref{sect_comp}, being a technical heart of the paper, we
prove, in a sequence of steps, that the lists of representatives
of the tilting-cotilting equivalence classes of gentle two-cycle
algebras given in Theorems~\ref{theo_nondeg} and~\ref{theo_deg}
are complete, while in Section~\ref{sect_min} we show that
different algebras from the list given in
Theorem~\ref{theo_nondeg} are not derived equivalent. The last
property follows from calculations of the derived invariant
introduced by Avella-Alaminos and Geiss in~\cite{AAGei2006}.

For a basic background on representation theory of finite
dimensional algebras we refer to~\cite{AssSimSko2006}.

Authors gratefully acknowledges the support from the Polish
Scientific Grant KBN No.~1 P03A 018 27.

\section{Basic tools and auxiliary results} \label{sect_prem}

By a (finite) quiver $\Delta$ we mean a finite set $\Delta_0$ of
vertices together with a finite set $\Delta_1$ of arrows and two
maps $s = s_\Delta, t = t_\Delta : \Delta_1 \to \Delta_0$ which
assign to an arrow $\alpha$ its starting and terminating vertex,
respectively. We say that an arrow $\alpha$ is adjacent to a
vertex $x$ if either $s \alpha = x$ or $t \alpha = x$. By a path
of length $n \in \bbN_+$ we mean a sequence $\sigma = \alpha_1
\cdots \alpha_n$ of arrows such that $s \alpha_i = t \alpha_{i +
1}$ for all $i \in [1, n - 1]$. In the above situation we denote
$s \alpha_n$ and $t \alpha_1$ by $s \sigma$ and $t \sigma$,
respectively. We also call $\alpha_1$ and $\alpha_n$ the
terminating and the starting arrow of $\sigma$, respectively.
Additionally, for each $x \in \Delta_0$ we consider the trivial
path of length $0$, also denoted by $x$, such that $s x = x = t
x$. The length of a path $\sigma$ will be denoted by $\ell
(\sigma)$. A path $\sigma$ is called maximal if there exists no
arrow $\alpha$ such that either $s \alpha = t \sigma$ or $t \alpha
= s \sigma$. Similarly, we define maximal paths starting (or
terminating) at a given vertex. A connected quiver is said to be
$c$-cycle if $|\Delta_1| = |\Delta_0| + c - 1$.

With a quiver $\Delta$ we associate its path algebra $k \Delta$,
which as a $k$-vector space has a basis formed by all paths in
$\Delta$ and whose multiplication is induced by the composition of
paths. By a relation $\rho$ in $\Delta$ we mean a linear
combination of paths of length at least $2$ with common starting
and terminating vertices. The common starting vertex is denoted by
$s \rho$ and the common terminating vertex by $t \rho$. A set $R$
of relations is called minimal if $\rho$ does not belong to the
ideal $\langle R \setminus \{ \rho \} \rangle$ of $k \Delta$
generated by $R \setminus \{ \rho \}$ for every $\rho \in R$. A
pair $(\Delta, R)$ consisting of a quiver $\Delta$ and a minimal
set of relations $R$ such that there exists $n \in \bbN$ with
$\sigma \in \langle R \rangle$ for each path $\sigma$ in $\Delta$
of length at least $n$, is called a bound quiver. If $(\Delta, R)$
is a bound quiver, then the algebra $k \Delta / \langle R \rangle$
is called the bound quiver algebra of $(\Delta, R)$.

Let $(\Delta, R)$ be a bound quiver and assume that $R$ consists
of paths. A path $\sigma$ in $\Delta$ is said to be a path
$(\Delta, R)$ if $\sigma \not \in \langle R \rangle$ (in other
words, none of the paths from $R$ is a subpath of $\sigma$). A
path $\sigma$ in $(\Delta, R)$ is said to be maximal if there is
no $\alpha \in \Delta_1$ such that either $s \alpha = t \sigma$
and $\alpha \sigma \not \in \langle R \rangle$ or $t \alpha = s
\sigma$ and $\sigma \alpha \not \in \langle R \rangle$. Again we
define maximal paths starting and terminating at a given vertex.
If additionally $R$ consists of paths of length two, then we say
that $\alpha \in \Delta_1$ is a free arrow provided there exists
no $\beta \in \Delta_1$ such that either $s \beta = t \alpha$ and
$\beta \alpha \in R$ or $t \beta = s \alpha$ and $\alpha \beta \in
R$.

Following~\cite{AssSko1987} we say that a connected bound quiver
$(\Delta, R)$ is gentle if the following conditions are satisfied:
\begin{enumerate}

\item
for each $x \in \Delta_0$ there are at most two arrows $\alpha$
such that $s \alpha = x$ ($t \alpha = x$),

\item
$R$ consists of paths of length two,

\item
for each $\alpha \in \Delta_1$ there is at most one arrow $\beta$
such that $t \beta = s \alpha$ and $\alpha \beta \not \in R$ ($s
\beta = t \alpha$ and $\beta \alpha \not \in R$),

\item
for each $\alpha \in \Delta_1$ there is at most one arrow $\beta$
such that $t \beta = s \alpha$ and $\alpha \beta \in R$ ($s \beta
= t \alpha$ and $\beta \alpha \in R$).

\end{enumerate}
An algebra which is isomorphic to the bound quiver algebra a
gentle bound quiver is called gentle.

With an abelian category $\calA$ we may associate its bounded
derived category $\calD^b (\calA)$ in the following way (see for
example~\cite{Ver1977} for details). The objects of $\calD^b
(\calA)$ are the bounded complexes of objects of $\calA$ and the
morphisms are obtained from the morphisms in the homotopy category
by formally inversing the quasi-isomorphisms (more precisely, by
localizing with respect to the quasi-isomorphism), where by a
quasi-isomorphism we mean a morphism of complexes which induces an
isomorphism of homology groups. The derived category together with
the shift functor sending $X$ to the shifted complex $X [1]$,
where $X [1]_n = X_{n + 1}$ and $d_{X [1]}^n = - d_{X [1]}^{n +
1}$ for $n \in \bbZ$, is a triangulated category (see for
example~\cite{Hap1988}). We say that two abelian categories
$\calA$ and $\calB$ are derived equivalent if there exists a
triangle equivalence $\calD^b (\calA) \to \calD^b (\calB)$. We say
that two algebras $\Lambda$ and $\Lambda'$ (bound quivers
$(\Delta, R)$ and $(\Delta', R')$) are derived equivalent if their
categories of modules (representations, respectively) are derived
equivalent. It follows from~\cite{SchZim2003}*{Corollary~1.2}
and~\cite{AAGei2006}*{Corollary~15} that for $c \in \bbZ$ gentle
$c$-cycle algebras (bound quivers) are closed under derived
equivalences.

Recall from~\cites{Bon1981, HapRin1982} that if $\Lambda$ is an
algebra, then a $\Lambda$-module $T$ is called tilting if
$\pd_\Lambda T \leq 1$, $\Ext_\Lambda^1 (T, T) = 0$, and $T$ is a
direct sum of $n$ pairwise nonisomorphic indecomposable modules,
where $n$ is the rank of the Grothendieck group of $\Lambda$.
Dually, we define notion of a cotilting module. Algebras $\Lambda$
and $\Lambda'$ are said to be tilting-cotilting equivalent if
there exists a sequence $\Lambda = \Lambda_0$, $\Lambda_1$,
\ldots, $\Lambda_n = \Lambda'$ of algebras such that for each $i
\in [0, n - 1]$ there exists a (co)tilting $\Lambda_{i -
1}$-module $T_{i - 1}$ such that $\Lambda_i \simeq
\End_{\Lambda_{i - 1}} (T_{i - 1})^{\op}$. It was proved by
Happel~\cite{Hap1987}*{Corollary~1.7} that if $\Lambda$ and
$\Lambda'$ are tilting-cotilting equivalent, then they are derived
equivalent.

A vertex $x$ in a quiver $\Delta$ is called a sink (source) if
there is no $\alpha \in \Delta_1$ with $s \alpha = x$ ($t \alpha =
x$, respectively). If $x$ is a sink in a gentle bound quiver
$(\Delta, R)$, then we define a new gentle bound quiver $(\Delta',
R')$, called the bound quiver obtained from $(\Delta, R)$ by
applying the APR-reflection at $x$, in the following way:
$\Delta'_0 = \Delta_0$, $\Delta'_1 = \Delta_1$,
\begin{align*}
s_{\Delta'} \alpha & =
\begin{cases}
x & \text{if } t_\Delta \alpha = x,
\\
s_\Delta \alpha & \text{otherwise},
\end{cases}
\\
t_{\Delta'} \alpha & =
\begin{cases}
s_\Delta \alpha & \text{if } t_\Delta \alpha = x,
\\
x & \text{if } \exists \; \beta \in \Delta_1 : t_\Delta \beta = x
\wedge s_\Delta \beta = t_\Delta \alpha \wedge \beta \alpha \in R,
\\
t_\Delta \alpha & \text{otherwise},
\end{cases}
\end{align*}
and
\begin{multline*}
R' = \{ \rho \in R \mid t_\Delta \rho \neq x \} \cup \{ \alpha
\beta \mid t_\Delta \alpha = x
\\
\wedge \exists \; \gamma \in \Delta_1 : \gamma \neq \alpha \wedge
t_\Delta \gamma = x \wedge s_\Delta \gamma = t_\Delta \beta \wedge
\gamma \beta \in R \}.
\end{multline*}
It follows that the bound quiver algebra of $(\Delta', R')$ is
isomorphic to the opposite algebra of the endomorphism algebra of
the APR-tilting module (see~\cite{AusPlaRei1979}) at $x$ defined
as
\[
\bigoplus_{\substack{a \in \Delta_0 \\ a \neq x}} P (a) \oplus
\Bigl( \bigoplus_{\substack{\alpha \in \Delta_1 \\ t \alpha = x}}
P (s \alpha) \Bigr) / P (x)
\]
(see~\cite{AssSko1987}*{2.1}).

We present now a generalization of the above construction due to
Brenner and Butler (see~\cite{BreBut1980}*{Chapter~2}). Let $x$ be
a vertex in a gentle bound quiver $(\Delta, R)$ such that there is
no $\alpha \in \Delta_1$ with $s \alpha = x = t \alpha$ and for
each $\alpha \in \Delta_1$ with $s \alpha = x$ there exists
$\beta_\alpha \in \Delta_1$ with $t \beta_\alpha = x$ and $\alpha
\beta_\alpha \not \in R$. We define a bound quiver $(\Delta', R')$
in the following way: $\Delta_0' = \Delta_0$, $\Delta_1' =
\Delta_1$,
\begin{align*}
s_{\Delta'} \alpha & =
\begin{cases}
x & \text{if } t_\Delta \alpha = x,
\\
s \beta_\alpha & \text{if } s_\Delta \alpha = x,
\\
s_\Delta \alpha & \text{otherwise},
\end{cases}
\\
t_{\Delta'} \alpha & =
\begin{cases}
s_\Delta \alpha & \text{if } t_\Delta \alpha = x,
\\
x & \text{if } \exists \; \beta \in \Delta_1 : t_\Delta \beta = x
\wedge s_\Delta \beta = t_\Delta \alpha \wedge \beta \alpha \in R,
\\
t_\Delta \alpha & \text{otherwise},
\end{cases}
\end{align*}
and set
\begin{multline*}
R' = \{ \alpha \beta \in R \mid t_\Delta \alpha \neq x \wedge
s_\Delta \alpha \neq x \} \cup \{ \alpha \beta_\alpha \mid
s_\Delta \alpha = x\} \cup
\\
\{ \alpha \beta \mid t_\Delta \alpha = x \wedge  \exists \; \gamma
\in \Delta_1 : \gamma \neq \alpha \wedge t_\Delta \gamma = x
\wedge s_\Delta \gamma = t_\Delta \beta \wedge \gamma \beta \in R
\}.
\end{multline*}
We will say that $(\Delta', R')$ is obtained from $(\Delta, R)$ by
applying the generalized APR-reflection at $x$. Similarly as in
the previous situation it follows easily that the bound quiver
algebra of $(\Delta', R')$ is the opposite algebra of the
endomorphism algebra of the tilting module defined in the same way
as before. Obviously all APR-reflections are examples of
generalized APR-reflections.

We also have a version of the above construction for a vertex $x$
of a gentle bound quiver $(\Delta, R)$ such that there exists
$\alpha \in \Delta_1$ with $s \alpha = x = t \alpha$. Observe that
then $\alpha^2 \in R$. We additionally assume that there exists
$\beta_0 \in \Delta_1$ with $s \beta_0 \neq x$ and $t \beta_0 =
x$. We define a bound quiver $(\Delta', R')$ in the following way:
$\Delta_0' = \Delta_0$, $\Delta_1' = \Delta_1$,
\begin{align*}
s_{\Delta'} \alpha & =
\begin{cases}
x & \text{if } t_\Delta \alpha = x,
\\
s_\Delta \beta_0 & \text{if } s_\Delta \alpha = x \wedge t_\Delta
\alpha \neq x,
\\
s_\Delta \alpha & \text{otherwise},
\end{cases}
\\
t_{\Delta'} \alpha & =
\begin{cases}
s_\Delta \alpha & \text{if } t_\Delta \alpha = x,
\\
x & \text{if } \exists \; \beta \in \Delta_1 : t_\Delta \beta = x
\wedge s_\Delta \beta = t_\Delta \alpha \wedge \beta \alpha \in R,
\\
t_\Delta \alpha & \text{otherwise},
\end{cases}
\end{align*}
and $R' = R$. We will say again that $(\Delta', R')$ is obtained
from $(\Delta, R)$ by applying the generalized APR-reflection at
$x$. It follows that the bound quiver algebra of $(\Delta', R')$
is the opposite algebra of the endomorphism algebra of the tilting
module
\[
\bigoplus_{\substack{a \in \Delta_0 \\ a \neq x}} P (a) \oplus (P
(y) \oplus P (y)) / P (x),
\]
where $y = s \beta_0$ and $P (x)$ is embedded in $P (y) \oplus P
(y)$ in such a way that the quotient module is indecomposable.

Let again $x$ be a sink in a gentle bound quiver $(\Delta, R)$. We
define the HW-reflection of $(\Delta, R)$ at $x$ as the bound
quiver $(\Delta', R')$ constructed in the following way. If
$\Delta_0 = \{ x \}$ (equivalently, $\Delta_1 = \varnothing$),
then $(\Delta', R') = (\Delta, R)$, hence assume this is not the
case. Then we put $\Delta_0' = \Delta_0$ and $\Delta_1' =
\Delta_1$. For each arrow $\alpha$ terminating at $x$ let
$\beta_\alpha$ be the starting arrow of the maximal path in
$(\Delta, R)$ terminating at $x$ whose terminating arrow is
$\alpha$. We put
\begin{align*}
s_{\Delta'} \alpha & =
\begin{cases}
x & \text{if } t_\Delta \alpha = x,
\\
s_\Delta \alpha & \text{otherwise},
\end{cases}
& \qquad \text{and} \qquad t_{\Delta'} \alpha & =
\begin{cases}
s_\Delta \beta_\alpha & \text{if } t_\Delta \alpha = x,
\\
t_\Delta \alpha & \text{otherwise}.
\end{cases}
\end{align*}
Finally let
\begin{multline*}
R' = \{ \rho \in R \mid t_\Delta \rho \neq x \}
\\
\cup \{ \beta \alpha \mid t_\Delta \alpha = x \wedge s_\Delta
\beta = s_\Delta \beta_\alpha \wedge \beta \neq \beta_\alpha
\wedge t_\Delta \beta \neq x \}.
\end{multline*}
It is known that bound quiver algebra of $(\Delta', R')$ is
(isomorphic to) the algebra obtained from the bound quiver algebra
of $(\Delta, R)$ by the HW-reflection at $x$ (defined
in~\cite{HugWas1983}), hence in particular is titling-cotilting
equivalent to $(\Delta, R)$ (see~\cite{TacWak1986}). Dually, one
defines a quiver obtained from $(\Delta, R)$ by applying the
HW-coreflection at a source.

Before we present basic applications of the above transformations,
we describe one more construction. Let $\Sigma$ be a subquiver of
a quiver $\Delta$. Assume that $\Sigma'$ is a quiver such that
$\Sigma_0' = \Sigma_0$ and $\Sigma_1' = \Sigma_1$ (but, usually,
$s_{\Sigma'} \neq s_\Sigma$ and $t_{\Sigma'} \neq t_\Sigma$). We
say that a quiver $\Delta'$ is obtained from $\Delta$ by replacing
$\Sigma$ by $\Sigma'$ if $\Delta'_0 = \Delta_0$, $\Delta'_1 =
\Delta_1$, and
\[
s_{\Delta'} \alpha =
\begin{cases}
s_\Delta \alpha & \text{if } \alpha \in \Delta_1 \setminus
\Sigma_1,
\\
s_{\Sigma'} \alpha & \text{if } \alpha \in \Sigma_1,
\end{cases}
\qquad t_{\Delta'} \alpha =
\begin{cases}
t_\Delta \alpha & \text{if } \alpha \in \Delta_1 \setminus
\Sigma_1,
\\
t_{\Sigma'} \alpha & \text{if } \alpha \in \Sigma_1,
\end{cases}
\]
for $\alpha \in \Delta_1$.

We describe now operations of shifting relations.

\begin{lemm} \label{lemm_shift_rel}
If
\[
\Sigma = \xymatrix{\vertexU{u} & \vertexU{x} \ar[l]_{\alpha_1} &
\vertexU{y} \ar[l]_{\alpha_2} & \vertexU{v} \ar[l]_{\alpha_3}}
\]
is a subquiver of a gentle bound quiver $(\Delta, R)$ such that
$\alpha_1 \alpha_2 \in R$, $\alpha_2 \alpha_3 \not \in R$, and
there are no other arrows adjacent to $y$, then $(\Delta, R)$ is
tilting-cotilting equivalent to the bound quiver $(\Delta', R')$,
where $R' = R \setminus \{ \alpha_1 \alpha_2 \} \cup \{ \alpha_2
\alpha_3 \}$ and $\Delta'$ is obtained from $\Delta$ by replacing
$\Sigma$ by the quiver
\[
\xymatrix{\vertexU{u} & \vertexU{y} \ar[l]_{\alpha_1} &
\vertexU{x} \ar[l]_{\alpha_2} & \vertexU{z} \ar[l]_{\alpha_3}}.
\]
\end{lemm}

\begin{proof}
Apply the generalized APR-coreflection at $y$.
\end{proof}

We remark that it may happen that one of the following equalities
hold: $u = y$, $x = v$ or $u = v$ holds. Moreover, $u = y$ if and
only if $x = v$, and in this case $\alpha_1 = \alpha_3$. We call
the above operation shifting the relation $\alpha_1 \alpha_2$ to
the right. Dually, one defines the operation of shifting relations
to the left.

We will need the following generalization of the above lemma.

\begin{lemm} \label{lemm_shift_last}
If
\[
\Sigma = \xymatrix{\vertexU{u} & \vertexU{x} \ar[l]_{\alpha_1} &
\vertexU{y_n} \ar[l]_{\alpha_2} \ar[r]^{\beta_n} & \vertexU{y_{n -
1}} \ar[rr]|{\textstyle \cdots} & & \vertexU{y_1} \ar[r]^{\beta_1}
& \vertexU{y_0} & \vertexU{v} \ar[l]_{\alpha_3}}, \; n \in \bbN_+,
\]
is a subquiver of a gentle bound quiver $(\Delta, R)$ such that
$\alpha_1 \alpha_2 \in R$, $\beta_1$, \ldots, $\beta_n$ are free
arrows, and there are no other arrows adjacent to $x$, $y_0$,
\ldots, $y_n$, then $(\Delta, R)$ is tilting-cotilting equivalent
to the bound quiver $(\Delta', R')$, where $R' = R \setminus \{
\alpha_1 \alpha_2 \} \cup \{ \alpha_2 \alpha_3 \}$ and $\Delta'$
is obtained from $\Delta$ by replacing $\Sigma$ by the quiver
\[
\xymatrix{\vertexU{u} & \vertexU{y_0} \ar[l]_{\alpha_1}
\ar[r]^{\beta_1} & \vertexU{y_1} \ar[rr]|{\textstyle \cdots} & &
\vertexU{y_{n - 1}} \ar[r]^{\beta_n} & \vertexU{y_n} & \vertexU{x}
\ar[l]_{\alpha_2} & \vertexU{z} \ar[l]_{\alpha_3}}.
\]
\end{lemm}

\begin{proof}
We leave it to the reader to verify that the following sequence of
operations leads from $(\Delta, R)$ to $(\Delta', R')$: first for
each $i = n, \ldots, 1$ we apply the APR-coreflections at $y_i$,
\ldots, $y_n$, $x$, and next we apply the generalized
APR-coreflections at $y_0$, \ldots, $y_n$.
\end{proof}

We will also shift a group of relations in the following sense.

\begin{lemm} \label{lemm_shift_rel_group}
Let
\[
\Sigma = \xymatrix{\vertexU{y} \ar[r]^\beta & \vertexU{x_0} &
\vertexU{x_1} \ar[l]_{\alpha_1} & & \vertexU{x_{n - 1}}
\ar[ll]|{\textstyle \cdots} & \vertexU{x_n} \ar[l]_{\alpha_n}}, \;
n \geq 2,
\]
be a subquiver of a gentle bound quiver $(\Delta, R)$ such that
$\beta$ is a free arrow, $\alpha_i \alpha_{i + 1} \in R$ for all
$i \in [1, n - 1]$, and there are no other arrows adjacent to
$x_0$, \ldots, $x_{n - 1}$. If there is no $\alpha \in \Delta_1$
with $t \alpha = x_n$ and $\alpha_n \alpha \in R$, then $(\Delta,
R)$ is tilting-cotilting equivalent to the bound quiver $(\Delta',
R')$, where $R' = R \setminus \{ \alpha_{n - 1} \alpha_n \} \cup
\{ \beta \alpha_1 \}$ and $\Delta'$ is obtained from $\Delta$ by
replacing $\Sigma$ by the quiver
\[
\xymatrix{\vertexU{y} & \vertexU{x_1} \ar[l]_\beta & \vertexU{x_2}
\ar[l]_{\alpha_1} & & \vertexU{x_{n - 1}} \ar[ll]|{\textstyle
\cdots} & \vertexU{x_0} \ar[l]_{\alpha_{n - 1}} \ar[r]^{\alpha_n}
& \vertexU{x_n}}.
\]
\end{lemm}

\begin{proof}
We apply the APR-reflection at $x_0$, followed by the composition
of the APR-reflection at $x_i$ and the generalized APR-reflection
at $x_0$ applied for $i = 1, \ldots, n - 1$.
\end{proof}

Observe that in the above lemma we shift relations to the left.
Dually we define an operation of shifting a group of relations to
the right.

We present now a reduction, being a consequence of the above
lemmas, which will appear a few times in our proofs. Let
\[
\Sigma = \xymatrix{\vertexU{x_0} \ar@{-}[r]^{\alpha_1} &
\vertexU{x_1} \ar@{-}[rr]|{\textstyle \cdots} & & \vertexU{x_{n -
1}} \ar@{-}[r]^{\alpha_n} & \vertexU{x_n}}, \; n \in \bbN_+,
\]
be a subquiver of a gentle bound quiver $(\Delta, R)$ such that
there are no other arrows adjacent to $x_1$, \ldots, $x_{n - 1}$
(it may happen that $x_0 = x_n$). We divide $\Sigma_1$ into two
disjoint subsets $\Sigma_{1, +}$ and $\Sigma_{1, -}$ in such a way
that, for each $i \in [1, n - 1]$, $\alpha_i$ and $\alpha_{i + 1}$
belong to the same subset if and only if either $s \alpha_i = t
\alpha_{i + 1}$ or $t \alpha_i = s \alpha_{i + 1}$. We
additionally assume that there exists $\varepsilon \in \{ -, + \}$
such that $\alpha \beta \not \in R$ for all $\alpha, \beta \in
\Sigma_{1, \varepsilon}$ with $s \alpha = t \beta$. If $x_0 = t
\alpha_1$, then by applying APR-reflections and shifts of
relations (we leave details to the reader), hence by passing to a
tilting-cotilting equivalent bound quiver, we may replace $\Sigma$
by the quiver
\[
\xymatrix{\vertexU{x_0} & \cdots \ar[l]_-{\alpha_1'} & \bullet
\ar[l]_-{\alpha_{l_1}'} \ar[r]^-{\alpha_{l_2}''} & \cdots
\ar[r]^-{\alpha_1''} & \bullet & \cdots \ar[l]_-{\alpha_1'''} &
\vertexU{x_n} \ar[l]_-{\alpha_{l_3}'''}}
\]
for some $l_1, l_2, l_3 \in \bbN$ with $l_1 + l_2 + l_3 = n$.
Moreover, we may additionally assume that $l_3 = 0$ if either $x_n
= t \alpha_n$ or $x_n = s \alpha_n$ and there is no $\alpha \in
\Delta_1$ with $t \alpha = x_n$ and $\alpha_n \alpha \in R$.
Obviously, we have the dual statement if $x_0 = s \alpha_1$.

The next observation is the following.

\begin{lemm} \label{lemm_op1}
If $p_1, p_2 \in \bbN_+$, $p_3, p_4 \in \bbN$, and $r_1 \in [0,
p_1 - 1]$, are such that $p_2 + p_3 \geq 2$ and $p_4 + r_1 \geq
1$, then $\Lambda_1 (p_1, p_2, p_3, p_4, r_1)$ and $\Lambda_1 (p_1
+ p_2 - r_1 - 1, r_1 + 1, p_4, p_3, p_2 - 1)$ are
tilting-cotilting equivalent.
\end{lemm}

\begin{proof}
It follows immediately by shifting relations.
\end{proof}

In order to formulate the next lemma we introduce a new family of
algebras. Namely, for $p_1, p_2 \in \bbN_+$, $p_3, p_4 \in \bbN$,
$r_1 \in [0, p_1 - 1]$, and $r_2 \in [0, p_2 - 1]$, such that $p_3
+ p_4 + r_1 + r_2 \geq 1$, let $\Lambda_2' (p_1, p_2, p_3, p_4,
r_1, r_2)$ be the algebra of the quiver
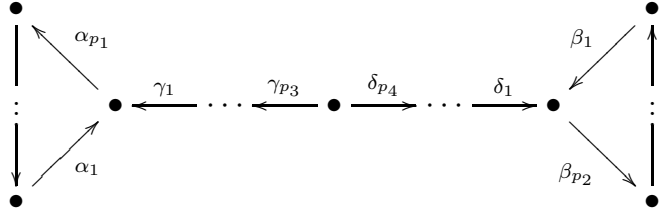
\[
\xymatrix{\bullet \ar[dd]|{\textstyle \vdots} & & & & & & \bullet
\ar[ld]_{\beta_1} \\ & \bullet \ar[lu]_{\alpha_{p_1}} & \cdots
\ar[l]_-{\gamma_1} & \bullet \ar[l]_-{\gamma_{p_3}}
\ar[r]^-{\delta_{p_4}} & \cdots \ar[r]^-{\delta_1} & \bullet
\ar[rd]_{\beta_{p_2}} \\ \bullet \ar[ru]_{\alpha_1} & & & & & &
\bullet \ar[uu]|{\textstyle \vdots}}
\]
bound by $\alpha_i \alpha_{i + 1}$ for $i \in [p_1 - r_1, p_1 -
1]$, $\alpha_{p_1} \alpha_1$, $\beta_i \beta_{i + 1}$ for $i \in
[p_2 - r_2, p_2 - 1]$, and $\beta_{p_2} \beta_1$.

\begin{lemm} \label{lemm_op2}
If $p_1, p_2, p_3 \in \bbN_+$, $p_4 \in \bbN$, $r_1 \in [0, p_1 -
1]$, and $r_2 \in [0, p_2 - 1]$, then $\Lambda_2' (p_1, p_2, p_3,
p_4, r_1, r_2)$ and $\Lambda_2' (p_1, p_2, p_3 - 1, p_4 + 1, r_1,
r_2)$ are tilting-cotilting equivalent.
\end{lemm}

\begin{proof}
Put $a_i = s \delta_i$, $i \in [1, p_4]$, and $b_i = s \beta_i$,
$i \in [1, p_2]$. We first apply the APR-coreflections at
$a_{p_4}$, \ldots, $a_1$, followed by the generalized
APR-coreflection at $b_{p_2}$ (we only apply the generalized
APR-coreflection at $b_{p_2}$ if $p_4 = 0$). Next we apply the
APR-coreflection at $b_{p_2 - i}$ followed by the generalized
APR-coreflection at $b_{p_2}$ for $i = 1, \ldots, r_2$ (we do
nothing in this step if $r_2 = 0$, hence in particular we do
nothing in this step if $p_2 = 1$), and finally we apply the
APR-coreflections at $b_{p_2 - r_2 - 1}$, \ldots, $b_1$ (there is
nothing to do if $r_2 = p_2 - 1$, hence again there is nothing to
do if $p_2 = 1$).
\end{proof}

\begin{coro} \label{coro_op2}
If $p_1, p_2 \in \bbN_+$, $p_3 \in \bbN$, $r_1 \in [0, p_1 - 1]$,
and $r_2 \in [0, p_2 - 1]$, are such that $p_3 + r_1 + r_2 \geq
1$, then $\Lambda_2 (p_1, p_2, p_3, r_1, r_2)$ and $\Lambda_2
(p_2, p_1, p_3, r_2, r_1)$ are tilting-cotilting equivalent.
\end{coro}

\begin{proof}
If follows immediately from the above lemma, since it is easily
seen that $\Lambda_2 (p_1, p_2, p_3, r_1, r_2)$ and $\Lambda_2
(p_2, p_1, p_3, r_2, r_1)$ are isomorphic to $\Lambda_2' (p_1,
p_2, p_3, 0, r_1, r_2)$ and $\Lambda_2' (p_1, p_2, 0, p_3, r_1,
r_2)$, respectively.
\end{proof}

\begin{prop} \label{prop_op}
If $\Lambda$ is one of the algebras listed in
Theorems~\ref{theo_nondeg} and~\ref{theo_deg}, then $\Lambda$ and
$\Lambda^{\op}$ are tilting-cotilting equivalent.
\end{prop}

\begin{proof}
If either $\Lambda = \Lambda_0 (p, r)$ for some $p \in \bbN_+$ and
$r \in [0, p - 1]$ or $\Lambda = \Lambda_0' (p, 0)$ for some $p
\in \bbN_+$, then the claim follows immediately by shifting
relations. If $\Lambda = \Lambda_1 (p_1, p_2, p_3, p_4, r_1)$ for
some $p_1, p_2 \in \bbN_+$, $p_3, p_4 \in \bbN$, and $r_1 \in [0,
p_1]$, such that $p_2 + p_3 \geq 2$ and $r_1 + p_4 \geq 1$, then
we have to additionally apply APR-coreflections. Finally, if
$\Lambda = \Lambda_2 (p_1, p_2, p_3, r_1, r_2)$ for some $p_1, p_2
\in \bbN_+$, $p_3 \in \bbN$, $r_1 \in [0, p_1 - 1]$, and $r_2 \in
[0, p_2 - 1]$, such that $p_3 + r_1 + r_2 \geq 1$, then the claim
follows either from Corollary~\ref{coro_op2}.
\end{proof}

An important consequence of the above lemma is that in our
considerations we may always replace an algebra by its opposite
algebra. Indeed, if for an algebra $\Gamma$ we are able to prove
that $\Gamma^{\op}$ is tilting-cotilting equivalent to an algebra
$\Lambda$ listed in Theorems~\ref{theo_nondeg} and~\ref{theo_deg},
then obviously $\Gamma$ is tilting-cotilting equivalent to
$\Lambda^{\op}$, hence also to $\Lambda$. In particular, once
Theorems~\ref{theo_nondeg} and~\ref{theo_deg} are proved, then we
know that if $\Gamma$ is a gentle two-cycle algebra, then $\Gamma$
and $\Gamma^{\op}$ are tilting-cotilting equivalent.

We finish this section by analyzing two particular families of
gentle two-cycle bound quivers. First, we prove the following.

\begin{prop} \label{prop_Lambda0}
If $(\Delta, R)$ is a gentle bound quiver such that
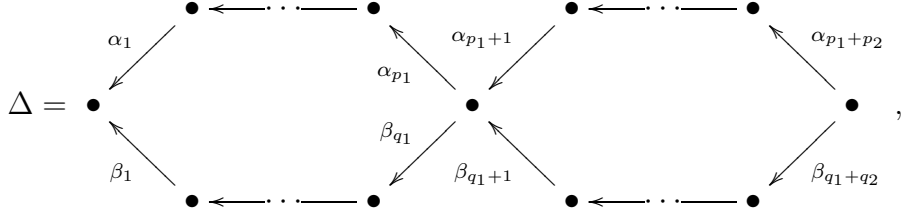
\[
\Delta = \vcenter{\xymatrix{& \bullet \ar[ld]_{\alpha_1} & &
\bullet \ar[ll]|{\textstyle \cdots} & & \bullet
\ar[ld]_{\alpha_{p_1 + 1}} & & \bullet \ar[ll]|{\textstyle \cdots}
\\ \bullet & & & & \bullet \ar[lu]^{\alpha_{p_1}}
\ar[ld]_{\beta_{q_1}} & & & & \bullet \ar[lu]_{\alpha_{p_1 + p_2}}
\ar[ld]^{\beta_{q_1 + q_2}} \\ & \bullet \ar[lu]^{\beta_1} & &
\bullet \ar[ll]|{\textstyle \cdots} & & \bullet
\ar[lu]^{\beta_{q_1 + 1}} & & \bullet \ar[ll]|{\textstyle
\cdots}}},
\]
for some $p_1, p_2, q_1, q_2 \in \bbN_+$, then the bound quiver
algebra of $(\Delta, R)$ is tilting-cotilting equivalent to
$\Lambda_0' (p, r)$ for some $p \in \bbN_+$ and $r \in [0, p -
1]$.
\end{prop}

We first show that also in the proof of this theorem we may pass
to the opposite algebras.

\begin{lemm} \label{lemm_extlist}
If $p \in \bbN_+$ and $r \in [1, p - 1]$, then $\Lambda_0' (p, r)$
and $\Lambda_0 (p + 1, r - 1)$ are tilting-cotilting equivalent.
\end{lemm}

\begin{proof}
In order to prove this equivalence, we put $x = t \beta$, $z = s
\delta$ and $x_1 = s \alpha_1$, and apply the APR-reflection at
$x$ followed by the APR-coreflection at $z$ and the APR-reflection
at $x_1$ to $\Lambda_0' (p, r)$. Then the claim follows by
shifting relations.
\end{proof}

\begin{coro}
If $p \in \bbN_+$ and $r \in [0, p - 1]$, then $\Lambda_0' (p, r)$
and $\Lambda_0' (p, r)^{\op}$ are tilting-cotilting equivalent.
\end{coro}

\begin{proof}
It follows either from Proposition~\ref{prop_op} (if $r = 0$) or
from the previous lemma and Proposition~\ref{prop_op} (if $r >
0$).
\end{proof}

In the proof of Proposition~\ref{prop_Lambda0} we will need the
following families of algebras:
\begin{itemize}

\item
$\Gamma_0(p, q, r)$ for $p, q \in \bbN_+$ and $r \in [0, p - 1]$
is the algebra of the quiver
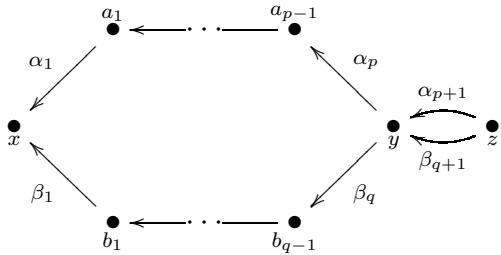
\[
\xymatrix{& \vertexD{a_1} \ar[ld]_{\alpha_1} & & \vertexD{a_{p -
1}} \ar[ll]|{\textstyle \cdots} \\ \vertexU{x} & & & & \vertexU{y}
\ar[lu]_{\alpha_p} \ar[ld]^{\beta_q} & \vertexU{z}
\ar@/_/[l]_{\alpha_{p + 1}} \ar@/^/[l]^{\beta_{q + 1}} \\ &
\vertexU{b_1} \ar[lu]^{\beta_1} & & \vertexU{b_{q - 1}}
\ar[ll]|{\textstyle \cdots}}
\]
bound by $\alpha_i \alpha_{i + 1}$ for $i \in [p - r, p]$ and
$\beta_q \beta_{q + 1}$,

\item
$\Gamma_1 (p, q, r, r')$ for $p, q \in \bbN_+$, $r \in [0, p -
1]$, and $r' \in \bbN$, is the algebra of the quiver
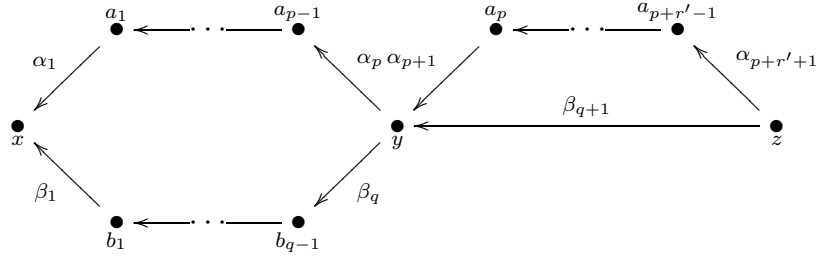
\[
\xymatrix{& \vertexD{a_1} \ar[ld]_{\alpha_1} & & \vertexD{a_{p -
1}} \ar[ll]|{\textstyle \cdots} & & \vertexD{a_p}
\ar[ld]_{\alpha_{p + 1}} & & \vertexD{a_{p + r' - 1}}
\ar[ll]|{\textstyle \cdots} \\ \vertexU{x} & & & & \vertexU{y}
\ar[lu]_{\alpha_p} \ar[ld]^{\beta_q} & & & & \vertexU{z}
\ar[lu]_{\alpha_{p + r' + 1}} \ar[llll]_{\beta_{q + 1}} \\ &
\vertexU{b_1} \ar[lu]^{\beta_1} & & \vertexU{b_{q - 1}}
\ar[ll]|{\textstyle \cdots}}
\]
bound by $\alpha_i \alpha_{i + 1}$ for $i \in [p - r, p + r']$ and
$\beta_q \beta_{q + 1}$,

\item
$\Gamma_2 (p, q, r, r')$ for $p, q \in \bbN_+$,  $r \in [0, p -
1]$, and $r' \in \bbN$, is the algebra of the quiver
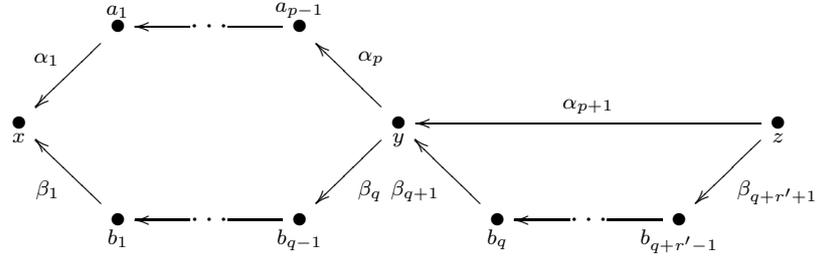
\[
\xymatrix{& \vertexD{a_1} \ar[ld]_{\alpha_1} & & \vertexD{a_{p -
1}} \ar[ll]|{\textstyle \cdots}
\\ \vertexU{x} & & & & \vertexU{y} \ar[lu]_{\alpha_p}
\ar[ld]^{\beta_q} & & & & \vertexU{z} \ar[ld]^{\beta_{q + r' + 1}}
\ar[llll]_{\alpha_{p + 1}} \\ & \vertexU{b_1} \ar[lu]^{\beta_1} &
& \vertexU{b_{q - 1}} \ar[ll]|{\textstyle \cdots} & &
\vertexU{b_q} \ar[lu]^{\beta_{q + 1}} & & \vertexU{b_{q + r' - 1}}
\ar[ll]|{\textstyle \cdots}}
\]
bound by $\alpha_i \alpha_{i + 1}$ for $i \in [p - r, p]$ and
$\beta_i \beta_{i + 1}$ for $i \in [q, q + r']$,

\end{itemize}
and the following series of lemmas.

\begin{lemm}
If $p, q \in \bbN_+$, $r \in [0, p - 1]$, and $q > 1$, then
$\Gamma_0(p, q, r)$ is tilting-cotilting equivalent to $\Gamma_0(p
+ 1, q - 1, r)$.
\end{lemm}

\begin{proof}
It is enough to apply the generalized APR-reflection at $b_{q -
1}$, followed by the APR-coreflection at $z$, the generalized
APR-coreflection at $y$, and the APR-coreflections at $b_{q - 2}$,
\ldots, $b_1$ (we omit the last step if $q = 2$).
\end{proof}

\begin{lemm}
If $p, q, \in \bbN_+$, $r \in [0, p - 1]$, $r' \in \bbN$, and $r'
\geq r$, then $\Gamma_1 (p, q, r, r')$ is tilting-cotilting
equivalent to $\Gamma_2 (q + r' - r, p, r' - r, r)$.
\end{lemm}

\begin{proof}
First for each $i \in [1, r]$ we apply the HW-coreflection at $z$
followed by the APR-reflection at $z$, and the generalized
APR-co\-ref\-lec\-tion at $a_{p + r' - i}$ applied $r + r' + 1 -
i$ times. Next we apply the HW-coreflections at $z$, $a_{p + r' -
r - 1}$, \ldots, $a_p$ (only at $z$ if $r = r'$) and we obtain a
bound quiver whose bound quiver algebra is easily seen to be
tilting-cotilting equivalent to $\Gamma_2 (q + r' - r, p, r' - r,
r)$.
\end{proof}

\begin{lemm}
If $p, q, \in \bbN_+$, $r \in [0, p - 1]$, $r' \in \bbN$, and $r
\geq r'$, then $\Gamma_1 (p, q, r, r')$ is tilting-cotilting
equivalent to $\Gamma_2 (p + 2 r' - r, q, r', r - r')$.
\end{lemm}

\begin{proof}
Since $\Gamma_1 (p, q, r, r')$ is tilting-cotilting equivalent to
$\Gamma_1 (p + r' - r, q, r', r)^{\op}$ and $\Gamma_2 (p + 2 r' -
r, q, r', r - r')$ is tilting-cotilting equivalent to $\Gamma_2 (q
+ r - r', p + r' - r, r - r', r')^{\op}$, hence the claim follows
from the previous lemma.
\end{proof}

\begin{lemm}
If $p, q \in \bbN_+$,  $r \in [0, p - 1]$, $r' \in \bbN$, and $r
\geq r'$, then $\Gamma_2 (p, q, r, r')$ is tilting-cotilting
equivalent to $\Gamma_2 (p, q, r - r', r')$.
\end{lemm}

\begin{proof}
By applying the APR-coreflection at $z$ followed by the
generalized APR-coreflection at $z$ applied $r'$ times, we replace
$\Gamma_2 (p, q, r, r')$ by (an algebra isomorphic to) the bound
quiver algebra of the quiver
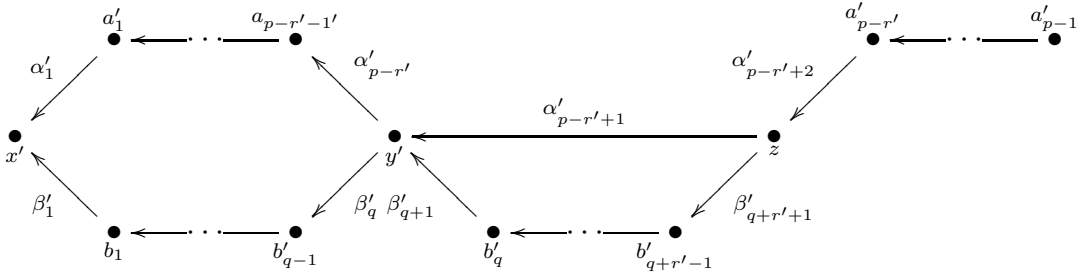
\[
\makebox[0pt]{\xymatrix{& \vertexD{a_1'} \ar[ld]_{\alpha_1'} & &
\vertexD{a_{p - r' - 1'}} \ar[ll]|{\textstyle \cdots} & & & & & &
\vertexD{a_{p - r'}'} \ar[ld]_{\alpha_{p - r' + 2}'} & &
\vertexD{a_{p - 1}'} \ar[ll]|{\textstyle \cdots} \\ \vertexU{x'} &
& & & \vertexU{y'} \ar[lu]_{\alpha_{p - r'}'} \ar[ld]^{\beta_q'} &
& & & \vertexU{z} \ar[ld]^{\beta_{q + r' + 1}'}
\ar[llll]_{\alpha_{p - r' + 1}'} \\ & \vertexU{b_1}
\ar[lu]^{\beta_1'} & & \vertexU{b_{q - 1}'} \ar[ll]|{\textstyle
\cdots} & & \vertexU{b_q'} \ar[lu]^{\beta_{q + 1}'} & &
\vertexU{b_{q + r' - 1}'} \ar[ll]|{\textstyle \cdots}}}
\]
bound by $\alpha_i' \alpha_{i + 1}'$ for $i \in [p - r, p]$ and
$\beta_i' \beta_{i + 1}'$ for $i \in [q, q + r']$. It is easily
seen that this algebra is tilting-cotilting equivalent to
$\Gamma_2 (p, q, r - r', r')$ (we just shift relations
sufficiently many times).
\end{proof}

\begin{lemm}
If $p, q \in \bbN$, $r \in [0, p - 1]$, $r' \in \bbN$, and $r'
\geq r$, then $\Gamma_2 (p, q, r, r')$ is tilting-cotilting
equivalent to $\Gamma_2 (p, q + r, r, r' - r)$.
\end{lemm}

\begin{proof}
Since $\Gamma_2 (p, q, r, r')$ is tilting-cotilting equivalent to
$\Gamma_2 (q + r', p - r, r', r)^{\op}$ and $\Gamma_2 (p, q + r,
r, r' - r)$ is tilting-cotilting equivalent to $\Gamma_2 (q + r',
p - r, r' - r, r)^{\op}$, hence the claim follows from the
previous lemma.
\end{proof}

\begin{proof}[Proof of Proposition~\ref{prop_Lambda0}]
Without loss of generality we may assume that $\alpha_{p_1}
\alpha_{p_1 + 1} \in R$ and $\beta_{q_1} \beta_{q_1 + 1} \in R$.
We first show that either $\alpha_i \alpha_{i + 1} \not \in R$ for
all $i \in [1, p_1 - 1]$ or $\beta_i \beta_{i + 1} \not \in R$ for
all $i \in [1, q_1 - 1]$. Assume this is not the case. In
particular, $p_1, q_1 \geq 2$. By shifting relations we may assume
that $\alpha_1 \alpha_2 \in R$ and $\beta_1 \beta_2 \in R$. If
$(\Delta', R')$ is the quiver obtained from $(\Delta, R)$ by
applying the HW-reflection at $x$ followed by the APR-reflection
at $x$, where $x = t \alpha_1$, then $\Delta' = \Delta$ and $R' =
R \setminus \{ \alpha_1 \alpha_2, \beta_1 \beta_2 \}$, hence the
claim follows by induction. Similarly, we prove that either
$\alpha_i \alpha_{i + 1} \not \in R$ for all $i \in [p_1 + 1, p_1
+ p_2 - 1]$ or $\beta_i \beta_{i + 1} \not \in R$ for all $i \in
[q_1 + 1, q_1 + q_2 - 1]$. Consequently, by shifting relations one
easily observes that the bound quiver algebra of $(\Delta, R)$ is
titling-cotilting equivalent either to $\Gamma_1 (p, q, r, r')$ or
to $\Gamma_2 (p, q, r, r')$ for some $p, q \in \bbN_+$, $r \in [0,
p - 1]$, and $r' \in \bbN$. Since $\Gamma_1 (p, q, r, 0) =
\Gamma_0 (p, q, r) = \Gamma_2 (p, q, r, 0)$ for all $p, q \in
\bbN_+$ and $r \in \bbN$, $\Gamma_1 (p, q, 0, r') \simeq \Gamma_0
(p + r', q, r')^{\op}$ and $\Gamma_2 (p, q, 0, r') \simeq \Gamma_0
(q + r', p, r')^{\op}$ for all $p, q \in \bbN_+$ and $r' \in
\bbN$, and $\Gamma_0 (p, 1, r)$ is tilting-cotilting equivalent to
$\Lambda_0' (p, r)$ for all $p \in \bbN_+$ and $r \in [0, p - 1]$,
hence the claim follows from the above series of lemmas.
\end{proof}

We finish this section with the following.

\begin{prop} \label{prop_Lambda0prim}
If $(\Delta, R)$ is a gentle bound quiver such that
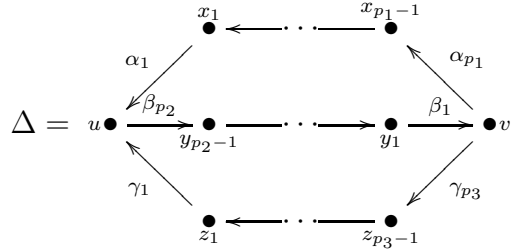
\[
\Delta = \vcenter{\xymatrix{& \vertexD{x_1} \ar[ld]_{\alpha_1} & &
\vertexD{x_{p_1 - 1}} \ar[ll]|{\textstyle \cdots} \\ \vertexR{u}
\ar[r]^{\beta_{p_2}} & \vertexU{y_{p_2 - 1}} \ar[rr]|{\textstyle
\cdots} & & \vertexU{y_1} \ar[r]^{\beta_1} & \vertexL{v}
\ar[lu]_{\alpha_{p_1}} \ar[ld]^{\gamma_{p_3}} \\ & \vertexU{z_1}
\ar[lu]^{\gamma_1} & & \vertexU{z_{p_3 - 1}} \ar[ll]|{\textstyle
\cdots}}}
\]
for some $p_1, p_2, p_3 \in \bbN_+$, and $\beta_{p_2} \alpha_1,
\gamma_{p_3} \beta_1 \in R$, then the bound quiver algebra of
$(\Delta, R)$ is tilting-cotilting equivalent to $\Lambda_0 (p,
r)$ for some $p \in \bbN_+$ and $r \in [0, p - 1]$.
\end{prop}

\begin{proof}
Let $r_1$ be the number of $i \in [1, p_1 - 1]$ such that
$\alpha_i \alpha_{i + 1} \in R$, let $r_2$ be the number of $i \in
[1, p_2 - 1]$ such that $\beta_i \beta_{i + 1} \in R$, and let
$r_3$ be the number of $i \in [1, p_3 - 1]$ such that $\gamma_i
\gamma_{i + 1} \in R$. We prove the claim by induction on $r_1 +
r_2 + r_3$.

If $r_1 = 0 = r_3$, then it follows by shifting relations that the
bound quiver algebra of $(\Delta, R)$ is tilting-cotilting
equivalent to $\Lambda_0 (p_1 + p_2 + p_3 - 2, r_2)$.

If $r_1 > 0$ and $r_3 = 0$, then by shifting relations we may
assume that $p_3 = 1$ and $\alpha_1 \alpha_2 \in R$. If $(\Delta',
R')$ is the bound quiver obtained from $(\Delta, R)$ by applying
the generalized APR-reflection at $u$ followed by the
APR-reflection at $x_1$, then $R' = R \setminus \{ \alpha_1
\alpha_2, \beta_{p_2} \alpha_1, \gamma_1 \beta_1 \} \cup \{
\gamma_1 \alpha_2, \beta_{p_2} \gamma_1 \}$ and
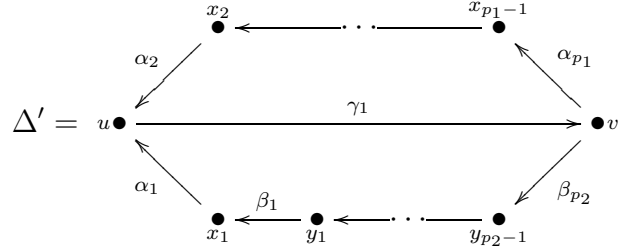
\[
\Delta' = \vcenter{\xymatrix{& \vertexD{x_2} \ar[ld]_{\alpha_2} &
& & \vertexD{x_{p_1 - 1}} \ar[lll]|{\textstyle \cdots} \\
\vertexR{u} \ar[rrrrr]^{\gamma_1} & & & & & \vertexL{v}
\ar[lu]_{\alpha_{p_1}} \ar[ld]^{\beta_{p_2}} \\ & \vertexU{x_1}
\ar[lu]^{\alpha_1} & \vertexU{y_1} \ar[l]_{\beta_1} & &
\vertexU{y_{p_2 - 1}} \ar[ll]|{\textstyle \cdots}}},
\]
hence the claim follows by induction. Dually, the claim follows if
$r_1 = 0$ and $r_3 > 0$.

Assume finally that $r_1 > 0$ and $r_3 > 0$. By shifting relations
we may assume that $\alpha_1 \alpha_2 \in R$ and $\gamma_1
\gamma_2 \in R$. If $(\Delta', R')$ is obtained from $(\Delta, R)$
by applying the generalized APR-reflection at $u$ followed by the
APR-reflection at $x_1$, then $R' = R \setminus \{ \alpha_1
\alpha_2, \beta_1 \alpha_1, \gamma_1 \gamma_2 \} \cup \{
\beta_{p_2} \gamma_1, \gamma_1 \alpha_2 \}$ and
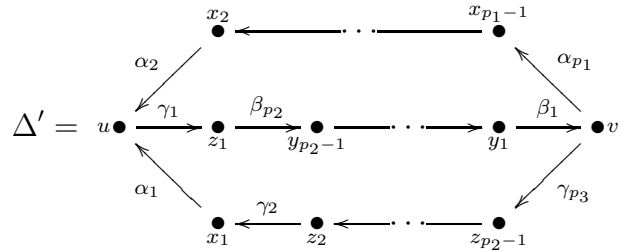
\[
\Delta' = \vcenter{\xymatrix{& \vertexD{x_2} \ar[ld]_{\alpha_2} &
& & \vertexD{x_{p_1 - 1}} \ar[lll]|{\textstyle \cdots} \\
\vertexR{u} \ar[r]^{\gamma_1} & \vertexU{z_1} \ar[r]^{\beta_{p_2}}
& \vertexU{y_{p_2 - 1}} \ar[rr]|{\textstyle \cdots} & &
\vertexU{y_1} \ar[r]^{\beta_1} & \vertexL{v}
\ar[lu]_{\alpha_{p_1}} \ar[ld]^{\gamma_{p_3}} \\ & \vertexU{x_1}
\ar[lu]^{\alpha_1} & \vertexU{z_2} \ar[l]_{\gamma_2} & &
\vertexU{z_{p_2 - 1}} \ar[ll]|{\textstyle \cdots}}},
\]
and the claim again follows by induction.
\end{proof}

\section{Completeness of the list} \label{sect_comp}

We start our considerations in this section by extending the list
of algebras in Theorems~\ref{theo_nondeg} and~\ref{theo_deg}.
Namely, as a consequence of the Lemmas~\ref{lemm_op1}
and~\ref{lemm_extlist} and Corollary~\ref{coro_op2}, it follows
that in order to show the completeness of the lists in
Theorems~\ref{theo_nondeg} and~\ref{theo_deg}, it is enough to
prove the following.

\begin{prop} \label{prop_mainprop}
If $(\Delta, R)$ is a gentle two-cycle bound quiver, then the
bound quiver algebra of $(\Delta, R)$ is tilting-cotilting
equivalent to one of the following algebras:
\begin{itemize}

\item
$\Lambda_0 (p, r)$ for some $p \in \bbN_+$ and $r \in [0, p - 1]$,

\item
$\Lambda_0' (p, r)$ for some $p \in \bbN_+$ and $r \in [0, p -
1]$,

\item
$\Lambda_1 (p_1, p_2, p_3, p_4, r_1)$ for some $p_1, p_2 \in
\bbN_+$, $p_3, p_4 \in \bbN$, and $r_1 \in [0, p_1 - 1]$, such
that $p_2 + p_3 \geq 2$ and $p_4 + r_1 \geq 1$,

\item
$\Lambda_2 (p_1, p_2, p_3, r_1, r_2)$ for some $p_1, p_2 \in
\bbN_+$, $p_3 \in \bbN$, $r_1 \in [0, p_1 - 1]$, $r_2 \in [0, p_2
- 1]$, such that $p_3 + r_1 + r_2 \geq 1$.

\end{itemize}
\end{prop}

For the rest of the section we assume that $(\Delta, R)$ is a
gentle two-cycle bound quiver. We show, in a sequence of steps,
that the bound quiver algebra of $(\Delta, R)$ is
tilting-cotilting equivalent to one of the algebras listed in the
above proposition.

We may divide the arrows in $\Delta$ into three disjoint groups:
\begin{itemize}

\item
$\alpha \in \Delta_1$ is called a cycle arrow if the quiver
$(\Delta_0, \Delta_1 \setminus \{ \alpha \})$ is connected,

\item
$\alpha \in \Delta_1$ is called a branch arrow if the quiver
$(\Delta_0, \Delta_1 \setminus \{ \alpha \})$ has a connected
component which is a two-cycle quiver,

\item
$\alpha \in \Delta_1$ is called a connecting arrow if the quiver
$(\Delta_0, \Delta_1 \setminus \{ \alpha \})$ has two connected
components which are one-cycle quivers.

\end{itemize}
A vertex $x$ of $\Delta$ is called a connecting vertex if there
exist at least three arrows adjacent to $x$ which are not branch
arrows. We call $\alpha \beta \in R$ a branch relation if either
$\alpha$ or $\beta$ is a branch arrow.

\begin{step} \label{step_one}
We may assume that there are no branch relations in $R$.
\end{step}

\begin{proof}
If there exists a branch relation in $(\Delta, R)$, then by
passing, if necessary, to the opposite algebra, we may assume that
there exists a subquiver
\[
\Sigma = \xymatrix{\vertexU{x_0} & \cdots \ar@{-}[l]_-{\alpha_1} &
\vertexU{x_{n - 2}} \ar@{-}[l]_-{\alpha_{n - 2}} & \vertexU{x_{n -
1}} \ar[l]_{\alpha_{n - 1}} & \vertexU{x_n} \ar[l]_{\alpha_n}},
\]
of $\Delta$ for some $n \geq 2$, where $\alpha_1$, \ldots,
$\alpha_{n - 2}$ are free arrows, $\alpha_{n - 1} \alpha_n \in R$,
and there are no other arrows adjacent to $x_0$, \ldots, $x_{n -
2}$ (in particular, $\alpha_{n - 1}$ is a branch arrow, hence
$\alpha_{n - 1} \alpha_n$ is a branch relation). By applying
APR-coreflections we may assume that $s \alpha_i = x_i$ for all $i
\in [1, n - 2]$. If $(\Delta', R')$ is the bound quiver obtained
from $(\Delta, R)$ by applying the generalized APR-reflections at
$x_{n - 2}$, \ldots, $x_1$ followed by the APR-reflection at
$x_0$, then $R' = R \setminus \{ \alpha_{n - 1} \alpha_n \}$ and
$\Delta'$ is obtained from $\Delta$ by replacing $\Sigma$ by the
quiver
\[
\xymatrix{\vertexU{x_{n - 1}} & \vertexU{x_0} \restore
\ar[l]_{\alpha_1} & \cdots \ar[l]_-{\alpha_2} & \vertexU{x_{n -
2}} \ar[l]_-{\alpha_{n - 1}} & \vertexU{x_n} \ar[l]_{\alpha_n}}.
\]
In particular, the number of branch relations decreases, hence the
claim follows by induction.
\end{proof}

By a branch in $\Delta$ we mean a maximal nontrivial (i.e.\ with
nonempty set of arrows) connected subquiver of $\Delta$ whose all
arrows are branch arrows. We say that a branch $B$ in $\Delta$ is
rooted at $x$ if $x \in B_0$ and there exists $\alpha \in
\Delta_1$ adjacent to $x$ which is not a branch arrow. An
immediate consequence of the assumption made in the above step is
that each branch $B$ in $\Delta$ is a linear quiver rooted at one
of its ends. Moreover, by applying the APR-reflections we may
assume that $B$ is equioriented and rooted at its unique sink.

\begin{step}
We may assume that there are no branch arrows in $\Delta$.
\end{step}

\begin{proof}
We say that $x \in \Delta_0$ is an insertion vertex if either $x$
is the connecting vertex or there exists $\alpha \in \Delta_1$
such that $s \alpha = x$, $\alpha$ is not a branch arrow, and
there is no $\beta \in \Delta_1$ with $t \beta = x$ and $\alpha
\beta \in R$. Observe that there is no branch which is rooted at
an insertion vertex. Moreover, for each $x \in \Delta_0$ there
exists a path in $\Delta$ starting at an insertion vertex and
terminating at $x$. In particular, if $B$ is a branch rooted at
$x$, then we call the minimal length of such a path the distance
between $B$ and an insertion vertex. We prove our claim by
induction on the number of branches in $(\Delta, R)$ and, for a
given branch $B$, by the distance between $B$ and an insertion
vertex.

Let
\[
B = \xymatrix{\vertexU{x_0} & \vertexU{x_1} \ar[l]_{\alpha_1} & &
\vertexU{x_{n - 1}} \ar[ll]|{\textstyle \cdots} & \vertexU{x_n}
\ar[l]_{\alpha_n}}, \; n \in \bbN_+,
\]
be a branch in $\Delta$. Let $\alpha$ and $\beta$ be the arrows in
$\Delta$ with $s \alpha = x_0 = t \beta$ and $\beta \neq
\alpha_1$. Observe that $\alpha \beta \in R$ and there are no
other arrows adjacent to $x_0$. Put $y = t \alpha$ and $z = s
\beta$.

Assume first that there is no $\gamma \in \Delta_1$ with $t \gamma
= z$ and $\beta \gamma \in R$. If $(\Delta', R')$ is the bound
quiver obtained from $(\Delta, R)$ by applying the generalized
APR-reflections at $x_0$, \ldots, $x_{n - 1}$, then $R' = R
\setminus \{ \alpha \beta \} \cup \{ \alpha \alpha_n \}$ and
$\Delta'$ is obtained from $\Delta$ by replacing the subquiver
\[
\xymatrix{& \vertexD{y} \\ \vertexU{z} \ar[r]^\beta &
\vertexU{x_0} \ar[u]_\alpha & \vertexU{x_1} \ar[l]_{\alpha_1} & &
\vertexU{x_{n - 1}} \ar[ll]|{\textstyle \cdots} & \vertexU{x_n}
\ar[l]_{\alpha_n}}
\]
by the quiver
\[
\xymatrix{\vertexU{z} & \vertexU{x_0} \ar[l]_\beta & \vertexU{x_1}
\ar[l]_{\alpha_1} & & \vertexU{x_{n - 1}} \ar[ll]|{\textstyle
\cdots}\ar[r]^{\alpha_n} & \vertexU{x_n} \ar[r]^\alpha &
\vertexU{y}},
\]
hence the claim follows in this case.

Assume now that there exists $\gamma \in \Delta_1$ with $t \gamma
= z$ and $\beta \gamma \in R$, and $z$ is a connecting vertex in
$\Delta_1$. Put $v = s \gamma$. If $(\Delta', R')$ is the bound
quiver obtained from $(\Delta, R)$ by applying the generalized
APR-reflections at $x_0$, \ldots, $x_{n - 1}$, then $R' = R
\setminus \{ \alpha \beta, \beta \gamma \} \cup \{ \alpha
\alpha_n, \alpha_n \gamma \}$ and $\Delta'$ is obtained from
$\Delta$ by replacing the subquiver
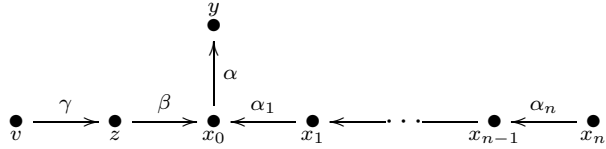
\[
\xymatrix{& & \vertexD{y} \\ \vertexU{v} \ar[r]^\gamma &
\vertexU{z} \ar[r]^\beta & \vertexU{x_0} \ar[u]_\alpha &
\vertexU{x_1} \ar[l]_{\alpha_1} & & \vertexU{x_{n - 1}}
\ar[ll]|{\textstyle \cdots} & \vertexU{x_n} \ar[l]_{\alpha_n}}
\]
by the quiver
\[
\vcenter
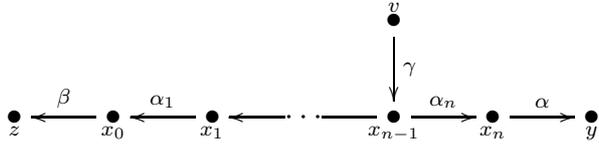
{\xymatrix{& & & & \vertexD{v} \ar[d]^\gamma \\
\vertexU{z} & \vertexU{x_0} \ar[l]_\beta & \vertexU{x_1}
\ar[l]_{\alpha_1} & & \vertexU{x_{n - 1}} \ar[ll]|{\textstyle
\cdots}\ar[r]^{\alpha_n} & \vertexU{x_n} \ar[r]^\alpha &
\vertexU{y}}}.
\]
Observe that the assumption that $z$ is a connecting vertex in
$\Delta$ implies that $\beta$, $\alpha_1$, \ldots, $\alpha_{n -
1}$ are not branch arrows in $\Delta'$.

Finally assume that there exists $\gamma \in \Delta_1$ with $t
\gamma = z$ and $\beta \gamma \in R$, but $z$ is not a connecting
vertex in $\Delta_1$. By induction we may assume that there is no
branch rooted at $z$. If $(\Delta', R')$ is the bound quiver
obtained from $(\Delta, R)$ by applying the HW-coreflection at
$x_i$ followed by the APR-reflection at $x_i$ for $i = n, \ldots,
1$, then $R' = R$ and $\Delta'$ is obtained from $\Delta$ by
replacing the subquiver
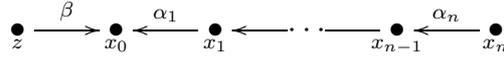
\[
\xymatrix{\vertexU{z} \ar[r]^\beta & \vertexU{x_0} & \vertexU{x_1}
\ar[l]_{\alpha_1} & & \vertexU{x_{n - 1}} \ar[ll]|{\textstyle
\cdots} & \vertexU{x_n} \ar[l]_{\alpha_n}}
\]
by the quiver
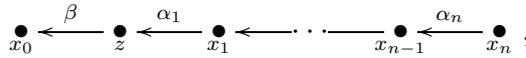
\[
\xymatrix{\vertexU{x_0} & \vertexU{z} \ar[l]_\beta & \vertexU{x_1}
\ar[l]_{\alpha_1} & & \vertexU{x_{n - 1}} \ar[ll]|{\textstyle
\cdots} & \vertexU{x_n} \ar[l]_{\alpha_n}},
\]
and the claim follows by induction.
\end{proof}

We say that $\Delta$ is special if either there is a unique
connecting vertex in $\Delta$ or there is a connecting arrow in
$\Delta$. Otherwise, we call $\Delta$ proper. We concentrate now
on the case when $\Delta$ is special. We describe first more
precisely its structure. We may divide the cycle arrows of
$\Delta$ into two disjoint subsets $\Delta_1^{(1)}$ and
$\Delta_1^{(2)}$ in such a way that cycle arrows $\alpha$ and $
\beta$ belong to the same subset if and only if the quiver
$(\Delta_0, \Delta_1 \setminus \{ \alpha, \beta \})$ has a
connected component which is a one-cycle quiver. For $j \in [1,
2]$ we denote by $\Delta^{(j)}$ the minimal subquiver of $\Delta$
with the set of arrows $\Delta_1^{(j)}$. Observe that
$\Delta^{(j)}$ is a (not necessarily oriented) cycle. We divide
the arrows in $\Delta^{(j)}$ into disjoint subsets $\Delta_{1,
-}^{(j)}$ and $\Delta_{1, +}^{(j)}$ in such a way that if $\alpha$
and $\beta$ are adjacent to the same vertex for $\alpha, \beta \in
\Delta_1^{(j)}$, $\alpha \neq \beta$, then $\alpha$ and $\beta$
belong to the same subset if and only if either $s \alpha = t
\beta$ or $t \alpha = s \beta$. For $\varepsilon \in \{ -, + \}$
we put
\[
R_\varepsilon^{(j)} = \{ \alpha \beta \in R \mid \alpha, \beta \in
\Delta_{1, \varepsilon}^{(j)} \}.
\]

\begin{step} \label{step_cycrel}
If $\Delta$ is special, then we may assume that for each $j \in
[1, 2]$ there exists $\varepsilon \in \{ -, + \}$ such that
$R_\varepsilon^{(j)} = \varnothing$.
\end{step}

\begin{proof}
If $\Delta^{(j)}$ is an oriented cycle, then there is nothing to
prove, hence assume that $\Delta^{(j)}$ is not an oriented cycle
and $R_-^{(j)} \neq \varnothing \neq R_+^{(j)}$. There exists a
subquiver
\[
\Sigma = \xymatrix{\vertexU{y_1} & \vertexU{y_2} \ar[l]_{\alpha_1}
& \vertexU{x_0} \ar[l]_{\alpha_2} & \cdots \ar@{-}[l]_-{\gamma_1}
\ar@{-}[r]^-{\gamma_n} & \vertexU{x_n} \ar[r]^{\beta_2} &
\vertexU{z_2} \ar[r]^{\beta_1} & \vertexU{z_1}}
\]
of $\Delta$ for some $n \in \bbN$, such that $\alpha_1 \alpha_2
\in R_-^{(j)}$, $\beta_1 \beta_2 \in R_+^{(j)}$, there are no
other arrows adjacent to $x_0$, \ldots, $x_n$, and $\gamma_1$,
\ldots, $\gamma_n$ are free arrows. By applying appropriate
APR-reflections at $x_1$, \ldots, $x_{n - 1}$ (see the discussion
after Lemma~\ref{lemm_shift_rel_group}) we may assume that
\[
\Sigma = \xymatrix{\vertexU{y_1} & \vertexU{y_2} \ar[l]_{\alpha_1}
& \vertexU{x_0} \ar[l]_{\alpha_2} & \cdots \ar[l]_-{\gamma_1} &
\vertexU{x_k} \ar[l]_-{\gamma_k} \ar[r]^-{\gamma_{k + 1}} & \cdots
\ar[r]^-{\gamma_n} & \vertexU{x_n} \ar[r]^{\beta_2} &
\vertexU{z_2} \ar[r]^{\beta_1} & \vertexU{z_1}}
\]
for some $k \in [0, n]$. By shifting the relations $\alpha_1
\alpha_2$ and $\beta_1 \beta_2$ to the right, we may assume that
$n = 0$, i.e.\
\[
\Sigma  = \xymatrix{\vertexU{y_1} & \vertexU{y_2}
\ar[l]_{\alpha_1} & \vertexU{x} \ar[l]_{\alpha_2} \ar[r]^{\beta_2}
& \vertexU{z_2} \ar[r]^{\beta_1} & \vertexU{z_1}}.
\]
Assume first that neither $y_2$ nor $z_2$ is a connecting vertex.
If $(\Delta', R')$ is the bound quiver obtained from $(\Delta, R)$
by applying the APR-co\-ref\-lec\-tions at $x$, $y_2$, and $z_2$,
then $R' = R \setminus \{ \alpha_1 \alpha_2, \beta_1 \beta_2 \}$
and $\Delta'$ is obtained from $\Delta$ by replacing $\Sigma$ by
the quiver
\[
\xymatrix{\vertexU{y_1} & \vertexU{z_2} \ar[l]_{\alpha_1} &
\vertexU{x} \ar[l]_{\beta_2} \ar[r]^{\alpha_2} & \vertexU{y_2}
\ar[r]^{\beta_1} & \vertexU{z_1}},
\]
and the claim follows by induction. Otherwise, we may assume
without loss of generality that $y_2$ is a connecting vertex and
$z_2$ is not a connecting vertex. If $(\Delta', R')$ is the bound
quiver obtained from $(\Delta, R)$ by applying the
APR-coreflections at $x$ and $z_2$, then $R' = R \setminus \{
\alpha_1 \alpha_2, \beta_1 \beta_2 \} \cup \{ \beta_1 \alpha_2 \}$
and $\Delta'$ is obtained from $\Delta$ by replacing $\Sigma$ by
the quiver
\[
\vcenter{\xymatrix{& & \vertexD{y_2} \ar[d]^{\alpha_2} \\
\vertexU{y_1} & \vertexU{z_2} \ar[l]_{\alpha_1} & \vertexU{x}
\ar[l]_{\beta_2} \ar[r]^{\beta_1} & \vertexU{z_1}}}.
\]
Observe that $\alpha_2$ is a connecting arrow in $\Delta'$, hence
the claim again follows by induction.
\end{proof}

\begin{step} \label{step_standard}
If $\Delta$ is special, then for each $j \in [1, 2]$ we may assume
that $\Delta^{(j)}$ is either an oriented cycle or there is a
unique source \textup{(}equivalently, unique sink\textup{)} in
$\Delta^{(j)}$.
\end{step}

\begin{proof}
This follows easily by applying APR-reflections and shifts of
relations (see the discussion after
Lemma~\ref{lemm_shift_rel_group}).
\end{proof}

\begin{step}
If $\Delta$ is special, then we may assume that either there is no
connecting arrow in $\Delta$ or, for each $j \in [1, 2]$,
$\Delta^{(j)}$ is an oriented cycle and $\alpha \beta \in R$ for
all $\alpha, \beta \in \Delta_1^{(j)}$ with $s \alpha = t \beta$.
\end{step}

\begin{proof}
We prove the claim by induction on the sum of the number of
connecting arrows and the number of connecting relations, where we
say that $\alpha \beta \in R$ is a connecting relation if both
$\alpha$ and $\beta$ are connecting arrows. We may assume without
loss of generality that either $\Delta^{(1)}$ is not an oriented
cycle or there exist $\alpha, \beta \in \Delta_1^{(1)}$ with $s
\alpha = t \beta$ and $\alpha \beta \not \in R$. Let $x \in
\Delta_0^{(1)}$ be a connecting vertex. Let $\alpha$ be the
connecting arrow adjacent to $x$. Without loss of generality we
may assume that $x = s \alpha$. Let $\beta$ and $\gamma$ be the
arrows adjacent to $x$ different from $\alpha$. Again we may
assume without loss of generality that $x = t \beta$. By symmetry
we may also assume that $\alpha \beta \in R$ if $x = t \gamma$.
Put $y = t \alpha$ and $z = s \beta$. In order to make it easier
to track the proof we will number the cases.

(1)~Assume that $\alpha \beta \not \in R$. According to our
assumptions this implies that $x = s \gamma$ and $\gamma \beta \in
R$. Put $v = t \gamma$. If $\Delta^{(1)}$ is not an oriented
cycle, then by applying APR-reflections and the dual of
Lemma~\ref{lemm_shift_last} we may assume that $v$ is a sink. In
particular, there is no $\gamma' \in \Delta_1$ with $s \gamma' =
v$ and $\gamma' \gamma \in R$. By shifting relations we may also
assume that this condition is satisfied if $\Delta^{(1)}$ is an
oriented cycle. Let $(\Delta', R')$ be the bound quiver obtained
from $(\Delta, R)$ by applying the generalized APR-coreflection at
$x$. If there is no $\alpha' \in \Delta_1$ with $s \alpha' = y$
and $\alpha' \alpha \in R$, then $R' = R \setminus \{ \gamma \beta
\} \cup \{ \alpha \beta \}$ and $\Delta'$ is obtained from
$\Delta$ by replacing the subquiver
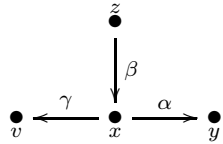
\[
\xymatrix{& \vertexD{z} \ar[d]^\beta \\ \vertexU{v} & \vertexU{x}
\ar[l]_\gamma \ar[r]^\alpha & \vertexU{y}}
\]
by the quiver
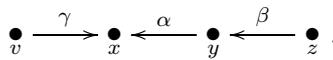
\[
\xymatrix{\vertexU{v} \ar[r]^\gamma & \vertexU{x} & \vertexU{y}
\ar[l]_\alpha & \vertexU{z} \ar[l]_\beta}.
\]
On the other hand, if there exists $\alpha' \in \Delta_1$ with $s
\alpha' = y$ and $\alpha' \alpha \in R$, then $R = R \setminus \{
\gamma \beta, \alpha' \alpha \} \cup \{ \alpha \beta, \alpha'
\gamma \}$ and $\Delta'$ is obtained from $\Delta$ by replacing
the subquiver
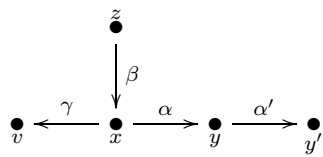
\[
\xymatrix{& \vertexD{z} \ar[d]^\beta \\ \vertexU{v} & \vertexU{x}
\ar[l]_\gamma \ar[r]^\alpha & \vertexU{y} \ar[r]^{\alpha'} &
\vertexU{y'}}
\]
by the quiver
\[
\vcenter
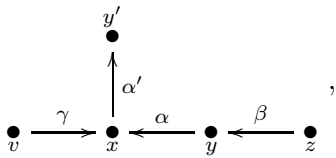
{\xymatrix{& \vertexD{y'} \\ \vertexU{v} \ar[r]^\gamma &
\vertexU{x} \ar[u]_{\alpha'} & \vertexU{y} \ar[l]_\alpha &
\vertexU{z} \ar[l]_\beta}},
\]
where $y' = t \alpha'$. Observe, that either $\Delta'$ is proper
(if $y$ is a connecting vertex in the second case) or we decrease
the number of connecting arrows (otherwise), hence the claim
follows by induction.

(2)~Assume that $\alpha \beta \in R$.

(2.1)~Assume that there is no $\alpha' \in \Delta_1'$ with $s
\alpha' = y$ and $\alpha' \alpha \in R$.

(2.1.1)~Assume that $y$ is a connecting vertex. If either
$\Delta^{(2)}$ is not an oriented cycle or there exist $\delta',
\delta'' \in \Delta_1^{(2)}$ with $s \delta' = t \delta''$ and
$\delta' \delta'' \not \in R$, then the claim follows by symmetry
from~(1), thus we may assume that $\Delta^{(2)}$ is an oriented
cycle such that $\delta' \delta'' \in R$ for all $\delta',
\delta'' \in \Delta_1^{(2)}$ with $s \delta' = t \delta''$.

(2.1.1.1)~Assume that $|\Delta_1^{(2)}| = 1$. If $(\Delta', R')$
is the bound quiver obtained from $(\Delta, R)$ by applying the
generalized APR-reflection at $y$, then $R' = R$ and $\Delta'$ is
obtained from $\Delta$ by replacing the subquiver
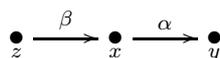
\[
\xymatrix{\vertexU{z} \ar[r]^\beta & \vertexU{x} \ar[r]^\alpha &
\vertexU{y}}
\]
by the quiver
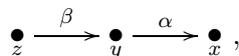
\[
\xymatrix{\vertexU{z} \ar[r]^\beta & \vertexU{y} \ar[r]^\alpha &
\vertexU{x}},
\]
hence the claim follows.

(2.1.1.2)~Assume that $|\Delta_1^{(2)}| > 1$. Let $\alpha'$ and
$\beta'$ be the arrows in $\Delta^{(2)}$ with $s \alpha' = y = t
\beta'$. Put $v' = t \alpha'$ and $x' = s \beta'$. Let $\gamma'$
be the arrow in $\Delta^{(2)}$ with $t \gamma' = x'$. Put $z' = s
\gamma'$. Recall that $\alpha' \beta', \beta' \gamma' \in R$. If
$(\Delta', R')$ is the bound quiver obtained from $(\Delta, R)$ by
applying the generalized APR-reflection at $y$ followed by the
APR-reflection at $x'$, then $R' = R \setminus \{ \alpha \beta,
\alpha' \beta', \beta' \gamma' \} \cup \{ \alpha \gamma', \alpha'
\alpha \}$ and $\Delta'$ is obtained from $\Delta$ by replacing
the subquiver
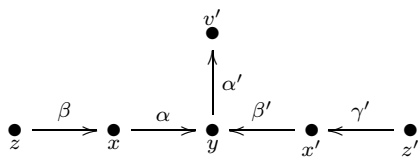
\[
\xymatrix{& & \vertexD{v'} \\ \vertexU{z} \ar[r]^\beta &
\vertexU{x} \ar[r]^\alpha & \vertexU{y} \ar[u]_{\alpha'} &
\vertexU{x'} \ar[l]_{\beta'} & \vertexU{z'} \ar[l]_{\gamma'}}
\]
by the quiver
\[
\vcenter{\xymatrix{& & \vertexD{z'} \ar[d]^{\gamma'} \\
\vertexU{z} \ar[r]^{\beta} & \vertexU{x'} \ar[r]^{\beta'} &
\vertexU{y} \ar[r]^\alpha & \vertexU{x} \ar[r]^{\alpha'} &
\vertexU{v'}}},
\]
hence the claim follows in this case.

(2.1.2)~Assume that $y$ is not a connecting vertex.

(2.1.2.1)~Assume that there exists $\alpha' \in \Delta_1$ with $s
\alpha' = y$. Our assumptions imply that $\alpha' \alpha \not \in
R$. Put $y' = t \alpha'$. If $(\Delta', R')$ is the bound quiver
obtained from $(\Delta, R)$ by applying the generalized
APR-reflection at $y$, then $R' = R \setminus \{ \alpha \beta \}
\cup \{ \alpha' \alpha \}$ and $\Delta'$ is obtained from $\Delta$
by replacing the subquiver
\[
\xymatrix{\vertexU{z} \ar[r]^\beta & \vertexU{x} \ar[r]^\alpha &
\vertexU{y} \ar[r]^{\alpha'} & \vertexU{y'}}
\]
by the quiver
\[
\xymatrix{\vertexU{z} \ar[r]^\beta & \vertexU{y} \ar[r]^\alpha &
\vertexU{x} \ar[r]^{\alpha'} & \vertexU{y'}},
\]
hence the claim follows by induction.

(2.1.2.2)~Assume there exists $\alpha' \in \Delta'$ with $t
\alpha' = y$. Put $x' = s \alpha'$.

(2.1.2.2.1)~Assume that either $x'$ is a connecting vertex or
$\alpha'$ is a free arrow. Moreover, if $x'$ is a connecting arrow
and $\alpha'$ is not a free arrow, then let $\beta'$ be the arrow
in $\Delta$ with $t \beta' = x'$ and $\alpha' \beta' \in R$, and
put $z' = s \beta'$. Let $(\Delta', R')$ be the bound quiver
obtained from $(\Delta, R)$ by applying the APR-reflection at $y$.
If $\alpha'$ is a free arrow, then $R' = R \setminus \{ \alpha
\beta \} \cup \{ \alpha' \beta \}$ and $\Delta'$ is obtained from
$\Delta$ by replacing the subquiver
\[
\xymatrix{\vertexU{z} \ar[r]^\beta & \vertexU{x} \ar[r]^\alpha &
\vertexU{y} & \vertexU{x'} \ar[l]_{\alpha'}}
\]
by the quiver
\[
\vcenter{\xymatrix{& \vertexD{z} \ar[d]^\beta \\ \vertexU{x} &
\vertexU{y} \ar[l]_\alpha \ar[r]^{\alpha'} & \vertexU{x'}}},
\]
hence the claim follows by induction. On the other hand, if $x'$
is a connecting arrow and $\alpha'$ is not a free arrow, then $R'
= R \setminus \{ \alpha \beta, \alpha' \beta' \} \cup \{ \alpha
\beta', \alpha' \beta \}$ and $\Delta'$ is obtained from $\Delta$
by replacing the subquiver
\[
\xymatrix{\vertexU{z} \ar[r]^\beta & \vertexU{x} \ar[r]^\alpha &
\vertexU{y} & \vertexU{x'} \ar[l]_{\alpha'} & \vertexU{z'}
\ar[l]_{\beta'}}
\]
by the quiver
\[
\vcenter{\xymatrix{\vertexD{z} \ar[rd]^\beta & & \vertexD{z'}
\ar[ld]_{\beta'} \\ \vertexU{x} & \vertexU{y} \ar[l]_\alpha
\ar[r]^{\alpha'} & \vertexU{x'}}},
\]
hence the claim follows.

(2.1.2.2.2)~Assume that $x'$ is not a connecting vertex and there
exists $\beta' \in \Delta_1$ with $t \beta' = x'$ and $\alpha'
\beta' \in R$. Put $z' = s \beta'$. If $(\Delta', R')$ is the
bound quiver obtained from $(\Delta, R)$ by applying the
APR-reflections at $y$ and $x'$, then $R' = R \setminus \{ \alpha
\beta, \alpha' \beta' \} \cup \{ \alpha \beta' \}$ and $\Delta'$
is obtained from $\Delta$ by replacing the subquiver
\[
\xymatrix{\vertexU{z} \ar[r]^\beta & \vertexU{x} \ar[r]^\alpha &
\vertexU{y} & \vertexU{x'} \ar[l]_{\alpha'} & \vertexU{z'}
\ar[l]_{\beta'}}
\]
by the quiver
\[
\vcenter{\xymatrix{& & \vertexD{x} \\ \vertexU{z} \ar[r]^\beta &
\vertexU{x'} \ar[r]^{\alpha'} & \vertexU{y} \ar[u]_\alpha &
\vertexU{z'} \ar[l]_{\beta'}}},
\]
hence the claim follows by induction.

(2.2)~Assume that there exists $\alpha' \in \Delta_1$ with $s
\alpha' = y$ and $\alpha' \alpha \in R$. Put $y' = t \alpha'$.

(2.2.1)~Assume that $x = t \gamma$. Let $\beta_1 \cdots \beta_n$
and $\gamma_1 \cdots \gamma_m$ be the maximal paths in $\Delta$
terminating at $x$ with $\beta_1 = \beta$ and $\gamma_1 = \gamma$.
Put $u = s \beta_n$, $u_i' = s \beta_i$ for $i \in [1, n - 1]$ and
$u_i'' = s \gamma_i$ for $i \in [1, m - 1]$.

(2.2.1.1)~Assume that there exists $i \in [1, m - 1]$ such that
$\gamma_i \gamma_{i + 1} \in R$. By shifting relations we may
assume that $\gamma_{m - 1} \gamma_m \in R$. Observe that $\beta_i
\beta_{i + 1} \not \in R$ for all $i \in [1, n - 1]$. If
$(\Delta', R')$ is the bound quiver obtained from $(\Delta, R)$ by
applying the HW-coreflection at $u$ followed by the composition of
the HW-coreflection at $u_i'$ and the APR-reflection at $u_i'$ for
$i = n - 1, \ldots, 1$, then $R' = R \setminus \{ \gamma_{m - 1}
\gamma_m, \alpha \beta \} \cup \{ \beta_n \gamma \}$ and $\Delta'$
is obtained from $\Delta$ by replacing the subquiver
\[
\xymatrix{\vertexU{u_{m - 1}''} & \vertexU{u} \ar[l]_{\gamma_m}
\ar[r]^{\beta_n} & \vertexU{u_{n  - 1}'} \ar[rr]|{\textstyle
\cdots} & & \vertexU{u_1'} \ar[r]^{\beta_1} & \vertexU{x}}
\]
by the quiver
\[
\xymatrix{\vertexU{u_{m - 1}''} \ar[r]^{\gamma_m} & \vertexU{u_{n
- 1}'} \ar[rr]|{\textstyle \cdots} & & \vertexU{u_1'}
\ar[r]^{\beta_1} & \vertexU{u} & \vertexU{x} \ar[l]_{\beta_n}},
\]
hence we reduce the proof to~(1).

(2.2.1.2)~Assume that $\gamma_i \gamma_{i + 1} \not \in R$ for all
$i \in [1, m - 1]$. Let $r$ be the number of $i \in [1, n - 1]$
such that $\beta_i \beta_{i + 1} \in R$. By shifting relations we
may assume $\beta_i \beta_{i + 1} \in R$ for all $i \in [n - r, n
- 1]$. Put $\beta_0 = \alpha$. If $(\Delta', R')$ is the bound
quiver obtained from $(\Delta, R)$ by applying the generalized
APR-coreflections at $u_1'$, \ldots, $u_{n - r - 1}'$, then $R' =
R \setminus \{ \alpha \beta \} \cup \{ \beta_{n - r - 1} \beta_{n
- r} \}$ and $\Delta'$ is obtained from $\Delta$ by replacing the
subquiver
\[
\xymatrix{\vertexU{u} \ar[r]^-{\beta_n} & \cdots \ar[r]^-{\beta_1}
& \vertexU{x} \ar[r]^{\beta_0} & \vertexU{y}}
\]
by the quiver
\[
\xymatrix{\vertexU{u} \ar[r]^-{\beta_n} & \cdots \ar[r]^-{\beta_{n
- r}} & \vertexU{x} \ar[r]^-{\beta_{n - r - 1}} & \cdots
\ar[r]^-{\beta_0} & \vertexU{y}}.
\]
Let $\gamma_1' \cdots \gamma_l'$ be the maximal path in $(\Delta',
R')$ with $\gamma_l' = \alpha$. Observe that $l > 1$ implies that
$y$ is a connecting vertex. Put
\[
u' =
\begin{cases}
u_{n - 1}' & r \geq 1,
\\
x & r = 0
\end{cases}
\qquad \text{and} \qquad v' = t \gamma_1'.
\]
Let $(\Delta'', R'')$ be the bound quiver obtained from $(\Delta',
R')$ by applying the HW-coreflection at $u$ followed by the
composition of the HW-coreflection at $u_i''$ and the
APR-reflection at $u_i''$ for $i = m - 1, \ldots, 1$. If there
exists $\delta$ in $\Delta$ with $t \delta = v'$ and $\delta \neq
\gamma_1'$, then $R'' = R' \setminus \{ \beta_{n - 1} \beta_n \}
\cup \{ \gamma_m \delta \}$, while $R'' = R' \setminus \{ \beta_{n
- 1} \beta_n \}$, otherwise. Moreover, $\Delta''$ is obtained from
$\Delta'$ by replacing the subquiver
\[
\xymatrix{\vertexU{u'} & \vertexU{u} \ar[l]_{\beta_n}
\ar[r]^-{\gamma_m} & \vertexU{u_{m - 1}''} \ar[r]^-{\gamma_{m -
1}} & \cdots \ar[r]^-{\gamma_1} & \vertexU{x} \ar[r]^-{\beta_{n -
r - 1}} & \cdots \ar[r]^-{\beta_0} & \vertexU{y}
\ar[r]^-{\gamma_{l - 1}'} & \cdots \ar[r]^{\gamma_1'} &
\vertexU{v'}}
\]
by the quiver
\[
\vcenter{\xymatrix{\vertexU{u'} \ar[r]^{\beta_n} & \vertexU{u_{m -
1}''} \ar[r]^-{\gamma_{m - 1}} & \cdots \ar[r]^{\gamma_1} &
\vertexU{u} & \vertexU{v'} \ar[l]_{\gamma_m} & \cdots
\ar[l]_-{\gamma_1'} & \vertexU{y} \ar[l]_-{\gamma_{l - 1}'} &
\cdots \ar[l]_-{\beta_0} & \vertexU{x} \ar[l]_-{\beta_{n - r -
1}}}},
\]
and the claim follows (by induction if $y$ is not a connecting
vertex).

(2.2.2)~Assume that $x = s \gamma$. Put $v = t \gamma$.

(2.2.2.1)~Assume that there exists $\gamma' \in \Delta_1$ with $s
\gamma' = v$ and $\gamma' \gamma \in R$ (by shifting relations we
may assume that this condition is satisfied if $\Delta^{(1)}$ is
an oriented cycle). Put $v' = t \gamma'$. If $(\Delta', R')$ is
the bound quiver obtained from $(\Delta, R)$ by applying the
generalized coreflection at $x$ followed, if $y$ is not a
connecting vertex, by the APR-coreflection at $y$, then
\[
R' =
\begin{cases}
R \setminus \{ \alpha' \alpha, \alpha \beta, \gamma' \gamma\} \cup
\{ \gamma' \alpha, \gamma \beta, \alpha' \gamma \} & \text{$y$ is
a connecting vertex},
\\
R \setminus \{ \alpha' \alpha, \alpha \beta, \gamma' \gamma\} \cup
\{ \gamma \beta, \alpha' \gamma \} & \text{\makebox[0pt][l]{$y$ is
not a connecting vertex,}}
\end{cases}
\]
and $\Delta'$ is obtained from $\Delta$ by replacing the subquiver
\[
\xymatrix{& & \vertexD{z} \ar[d]^\beta \\ \vertexU{v'} &
\vertexU{v} \ar[l]_{\gamma'} & \vertexU{x} \ar[l]_\gamma
\ar[r]^\alpha & \vertexU{y} \ar[r]^{\alpha'} & \vertexU{y'}}
\]
by the quiver
\[
\xymatrix{& \vertexD{v'} & & \vertexU{y'} \\ \vertexU{z}
\ar[r]^\beta & \vertexU{v} \ar[r]^\gamma & \vertexU{x}
\ar[ru]^{\alpha'} \ar[lu]_{\gamma'} & \vertexU{y} \ar[l]_\alpha}
\]
if $y$ is a connecting vertex, or
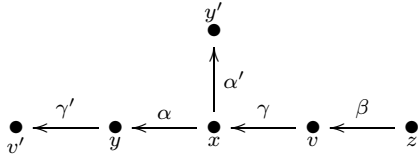
\[
\xymatrix{& & \vertexD{y'} \\ \vertexU{v'} & \vertexU{y}
\ar[l]_{\gamma'} & \vertexU{x} \ar[l]_\alpha \ar[u]_{\alpha'} &
\vertexU{v} \ar[l]_\gamma & \vertexU{z} \ar[l]_\beta}
\]
if $y$ is not a connecting vertex, hence the claim again follows.

(2.2.2.2)~Assume that $\Delta^{(1)}$ is not an oriented cycle. Let
$\gamma_1 \cdots \gamma_n$ be the maximal path in $\Delta$ with
$\gamma_n = \gamma$. We may additionally assume that $\gamma_i
\gamma_{i + 1} \not \in R$ for all $i \in [1, n - 1]$.
Consequently, we may reduce the proof in this case to~(2.2.1) by
applying APR-reflections and shifts of relations.
\end{proof}

\begin{step}
If $\Delta$ is special, then we may assume that for each $j \in
[1, 2]$ $\Delta^{(j)}$ is an oriented cycle or either the source
or the sink in $\Delta^{(j)}$ is a connecting vertex.
\end{step}

\begin{proof}
If both $\Delta^{(1)}$ and $\Delta^{(2)}$ are not oriented cycles,
then there is nothing to prove, thus without loss of generality we
may assume that $\Delta^{(1)}$ is not an oriented cycle. Observe
that our assumptions imply that there are no connecting arrows in
$\Delta$. Let $x$ be the connecting vertex in $\Delta$ and assume
that $x$ is neither a source nor a sink in $\Delta^{(1)}$. Observe
that $x \in \Delta_0^{(1)} \cap \Delta_0^{(2)}$. Let $\alpha$,
$\beta$, $\alpha'$ and $\beta'$ be the arrows in $\Delta$ with $s
\alpha = t \beta = x = s \alpha' = t \beta'$, $\alpha, \beta \in
\Delta_1^{(1)}$, and $\alpha', \beta' \in \Delta_1^{(2)}$. Put $y
= t \alpha$, $y' = t \alpha'$, $z = s \beta$, and $z' = s \beta'$.
By applying APR-coreflections, shifts of relations and
Lemma~\ref{lemm_shift_last} we may assume that $z$ is a source in
$\Delta^{(1)}$.

Assume first that $\alpha' = \beta'$. Then $\alpha \beta \in R$
and $\alpha' \beta' \in R$. Let $\gamma_1 \cdots \gamma_m$ be the
maximal path in $\Delta$ starting at $z$ with $\gamma_m \neq
\beta$. Observe that $\gamma_i \gamma_{i + 1} \not \in R$ for all
$i \in [1, m - 1]$. Put $v_i = s \gamma_i$ for $i \in [1, m - 1]$.
The bound quiver algebra of the bound quiver obtained from
$(\Delta, R)$ by applying the APR-coreflections at $z$, $v_{m -
1}$, \ldots, $v_1$, is easily seen to be tilting-cotilting
equivalent to $\Lambda_2 (p, 1, m, r, 0)$ for some $p \in \bbN_+$
and $r \in [0, p -1]$, hence the claim follows in this case.

Assume now that $\alpha \beta \in R$ and $\alpha' \beta' \in R$,
but $\alpha' \neq \beta'$. Let $(\Delta', R')$ be the bound quiver
obtained from $(\Delta, R)$ by applying the generalized
APR-reflection at $x$. If there exists $\beta'' \in \Delta_1$ with
$t \beta'' = z'$ and $\beta' \beta'' \in R$, then $R' = R
\setminus \{ \alpha \beta, \alpha' \beta', \beta' \beta'' \} \cup
\{ \alpha \beta', \alpha' \beta, \beta \beta'' \}$ and $\Delta'$
is obtained from $\Delta$ by replacing the subquiver
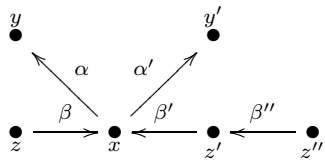
\[
\xymatrix{\vertexD{y} & & \vertexD{y'} \\ \vertexU{z} \ar[r]^\beta
& \vertexU{x} \ar[lu]_\alpha \ar[ru]^{\alpha'} & \vertexU{z'}
\ar[l]_{\beta'} & \vertexU{z''} \ar[l]_{\beta''}}
\]
by the quiver
\[
\vcenter{\xymatrix{& & \vertexD{z''} \ar[d]^{\beta''} \\
\vertexU{y} & \vertexU{z'} \ar[l]_\alpha & \vertexU{x}
\ar[l]_{\beta'} \ar[r]^\beta & \vertexU{z} \ar[r]^{\alpha'} &
\vertexU{y'}}},
\]
where $z'' = s \beta''$. Otherwise, $R' = R \setminus \{ \alpha
\beta, \alpha' \beta' \} \cup \{ \alpha \beta', \alpha' \beta \}$
and $\Delta'$ is obtained from $\Delta$ by replacing the subquiver
\[
\xymatrix{\vertexD{y} & & \vertexD{y'} \\ \vertexU{z} \ar[r]^\beta
& \vertexU{x} \ar[lu]_\alpha \ar[ru]^{\alpha'} & \vertexU{z'}
\ar[l]_{\beta'}}
\]
by the quiver
\[
\xymatrix{\vertexU{y} & \vertexU{z'} \ar[l]_\alpha & \vertexU{x}
\ar[l]_{\beta'} \ar[r]^\beta & \vertexU{z} \ar[r]^{\alpha'} &
\vertexU{y'}}.
\]
In particular, in both cases $\Delta'$ proper.

Assume finally that $\alpha \beta' \in R$ and $\alpha' \beta \in
R$. Let $\gamma$ be the arrow in $\Delta$ with $s \gamma = z$ and
$\gamma \neq \beta$. Put $v = t \gamma$. If there exists $\gamma'
\in \Delta_1$ with $s \gamma' = v$ and $\gamma' \gamma \in R$,
then let $(\Delta', R')$ be the bound quiver obtained from
$(\Delta, R)$ by applying the APR-coreflections at $z$ and $v$,
and let $v' = t \gamma'$. Observe that $R' = R \setminus \{
\alpha' \beta, \gamma' \gamma \} \cup \{ \gamma' \beta \}$ and
$\Delta'$ is obtained from $\Delta$ by replacing the subquiver
\[
\xymatrix{\vertexU{v'} & \vertexU{v} \ar[l]_{\gamma'} &
\vertexU{z} \ar[l]_\gamma \ar[r]^\beta & \vertexU{x}
\ar[r]^{\alpha'} & \vertexU{y'}}
\]
by the quiver
\[
\vcenter{\xymatrix{& \vertexD{x} \ar[d]^\beta \\ \vertexU{v'} &
\vertexU{z} \ar[l]_{\gamma'} \ar[r]^\gamma & \vertexU{v}
\ar[r]^{\alpha'} & \vertexU{y'}}}.
\]
Otherwise, if $(\Delta', R')$ is the bound quiver obtained from
$(\Delta, R)$ by applying the APR-coreflection at $z$, then $R' =
R \setminus \{ \alpha' \beta \} \cup \{ \alpha' \gamma \}$ and
$\Delta'$ is obtained from $\Delta$ by replacing the subquiver
\[
\xymatrix{\vertexU{v} & \vertexU{z} \ar[l]_\gamma \ar[r]^\beta &
\vertexU{x} \ar[r]^{\alpha'} & \vertexU{y'}}
\]
by the quiver
\[
\vcenter{\xymatrix{& \vertexD{y'} \\ \vertexU{v} \ar[r]^\gamma &
\vertexU{z} \ar[u]_{\alpha'} & \vertexU{x} \ar[l]_\beta}}.
\]
Again in both cases $\Delta'$ is proper and this finishes the
proof.
\end{proof}

\begin{step}
We may assume that $(\Delta, R)$ is proper.
\end{step}

\begin{proof}
If $\Delta^{(1)}$ is not an oriented cycle, then also
$\Delta^{(2)}$ is not an oriented cycle and the claim follows from
Proposition~\ref{prop_Lambda0}, thus assume that $\Delta^{(1)}$
(and consequently also $\Delta^{(2)}$) is an oriented cycle.

Assume first that there are no connecting arrows in $\Delta$ and
let $x$ be the connecting vertex in $\Delta$. Let $\alpha$,
$\beta$, $\alpha'$ and $\beta'$ be the arrows in $\Delta$ with $s
\alpha = t \beta = x = s \alpha' = t \beta'$, $\alpha, \beta \in
\Delta_1^{(1)}$, and $\alpha', \beta' \in \Delta_1^{(2)}$. If
$\alpha \beta \in R$ and $\alpha' \beta' \in R$, then it follows
by shifting relations that the bound quiver algebra of $(\Delta,
R)$ is tilting-cotilting equivalent to $\Lambda_2 (p_1, p_2, 0,
r_1, r_2)$ for some $p_1, p_2 \in \bbN_+$, $r_1 \in [0, p_1 - 1]$,
and $r_2 \in [0, p_2 - 1]$, such that $r_1 + r_2 \geq 1$. On the
other hand, if $\alpha \beta' \in R$ and $\alpha' \beta \in R$,
then it follows by shifting relations that the bound quiver
algebra of $(\Delta, R)$ is tilting-cotilting equivalent to
$\Lambda_1 (p_1, p_2, 0, 0, r_1)$ for some $p_1, p_2 \in \bbN_+$,
$p_1, p_2 \geq 2$, and $r_1 \in [1, p_1 - 1]$.

Now assume that there are connecting arrows in $\Delta$. Recall
that in this case $\alpha \beta \in R$ for all cycle arrows
$\alpha$ and $\beta$ with $s \alpha = t \beta$. Let $\Delta^{(0)}$
be the minimal subquiver of $\Delta$ with the set of arrows
consisting of the connecting arrows. Let $x \in \Delta_0^{(1)}$
and $y \in \Delta_0^{(2)}$ be the connecting vertices. Observe
that $\Delta^{(0)}$ is a linear quiver. We show that we may assume
that $x$ is a unique sink in $\Delta^{(0)}$, $y$ is a unique
source in $\Delta^{(0)}$, and there are no $\alpha, \beta \in
\Delta_1^{(0)}$ with $s \alpha = t \beta$ and $\alpha \beta \in
R$. This will immediately imply that the bound quiver algebra of
$(\Delta, R)$ is tilting-cotilting equivalent to $\Lambda_2 (p_1,
p_2, p_3, p_1 - 1, p_2 - 1)$ for some $p_1, p_2, p_3 \in \bbN_+$.

By repeating arguments from the proofs of Steps~\ref{step_cycrel}
and~\ref{step_standard} and passing, if necessary, to the opposite
algebra, we may assume that
\[
\Delta = \xymatrix{\vertexU{x} & \vertexU{x_1} \ar[l]_{\alpha_1} &
& \vertexU{x_{n - 1}} \ar[ll]|{\textstyle \cdots} & \vertexU{z}
\ar[l]_{\alpha_n} \ar[r]^{\beta_m} & \vertexU{y_{m - 1}}
\ar[rr]|{\textstyle \cdots} & & \vertexU{y_1} \ar[r]^{\beta_1} &
\vertexU{y}}
\]
for some $n \in \bbN_+$ and $m \in \bbN$, and $\beta_i \beta_{i +
1} \not \in R$ for all $i \in [1, m - 1]$. It is enough to show
that we may additionally assume that $\alpha_i \alpha_{i + 1} \not
\in R$ for all $i \in [1, n - 1]$, since then the claim follows
from Lemma~\ref{lemm_op2}. Assume this is not the case. By
shifting relations we may assume that $\alpha_1 \alpha_2 \in R$.

If $|\Delta^{(1)}_1| = 1$ and $(\Delta', R')$ is the bound quiver
obtained from $(\Delta, R)$ by applying the generalized
APR-reflection at $x$ followed by the APR-reflection at $x_1$,
then $R' = R \setminus \{ \alpha_1 \alpha_2 \}$ and $\Delta' =
\Delta$. Otherwise, let $\gamma$, $\delta$ and $\delta'$ be the
arrows in $\Delta^{(1)}$ with $s \gamma = x = t \delta$ and $t
\delta' = s \delta$. Observe that our assumptions imply that
$\gamma \delta, \delta \delta' \in R$. Put $u = t \gamma$, $v = s
\delta$ and $v' = s \delta'$. If $(\Delta', R')$ is the bound
quiver obtained from $(\Delta, R)$ by applying the generalized
APR-reflection at $x$ followed by the APR-reflection at $v$, then
$R' = R \setminus \{ \gamma \delta, \delta \delta', \alpha_1
\alpha_2 \} \cup \{ \gamma \alpha_1, \alpha_1 \delta' \}$ and
$\Delta'$ is obtained from $\Delta$ by replacing the subquiver
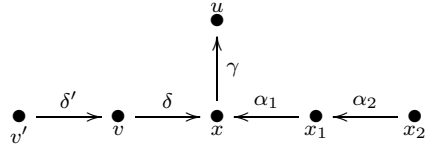
\[
\xymatrix{& & \vertexD{u} \\ \vertexU{v'} \ar[r]^{\delta'} &
\vertexU{v} \ar[r]^\delta & \vertexU{x} \ar[u]_\gamma &
\vertexU{x_1} \ar[l]_{\alpha_1} & \vertexU{x_2} \ar[l]_{\alpha_2}}
\]
by the quiver
\[
\vcenter
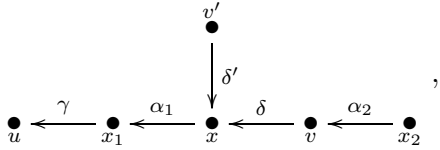
{\xymatrix{& & \vertexD{v'} \ar[d]^{\delta'} \\
\vertexU{u} & \vertexU{x_1} \ar[l]_\gamma & \vertexU{x}
\ar[l]_{\alpha_1} & \vertexU{v} \ar[l]_\delta & \vertexU{x_2}
\ar[l]_{\alpha_2}}},
\]
where $x_2 = z$ if $n = 2$. Consequently, in both cases the claim
follows by induction.
\end{proof}

We investigate now the case when $\Delta$ is proper. In this case
we may divide the arrows in $\Delta$ into there disjoint subsets
$\Delta_1^{(1)}$, $\Delta_1^{(2)}$, $\Delta_1^{(3)}$ in such a way
that $\alpha, \beta \in \Delta_1$ belong to the same subset if and
only if the quiver $(\Delta_0, \Delta_1 \setminus \{ \alpha, \beta
\})$ has a connected component which is a one-cycle quiver. For $j
\in [1, 3]$ we denote by $\Delta^{(j)}$ the minimal subquiver of
$\Delta$ with the set of arrows $\Delta_1^{(j)}$. Observe that
$\Delta^{(j)}$ is a linear quiver. We divide the arrows in
$\Delta^{(j)}$ into disjoint subsets $\Delta_{1, -}^{(j)}$ and
$\Delta_{1, +}^{(j)}$ in such a way that if $\alpha$ and $\beta$
are adjacent to the same vertex for $\alpha, \beta \in
\Delta_1^{(j)}$, $\alpha \neq \beta$, then $\alpha$ and $\beta$
belong to the same subset if and only if either $s \alpha = t
\beta$ or $t \alpha = s \beta$. For $j \in [1, 3]$ and
$\varepsilon \in \{ -, + \}$ we put
\[
R_\varepsilon^{(j)} = \{ \alpha \beta \in R \mid \alpha, \beta \in
\Delta_{1, \varepsilon}^{(j)} \}.
\]

\begin{step}
We may assume that either $R_+^{(j)} = \varnothing$ or $R_-^{(j)}
= \varnothing$ for each $j \in [1, 3]$.
\end{step}

\begin{proof}
Analogous to the proof of Step~\ref{step_cycrel}.
\end{proof}

\begin{step}
We may assume that either there is at most one sink in
$\Delta^{(j)}$ or there is at most one source in $\Delta^{(j)}$
for each $j \in [1, 3]$
\end{step}

\begin{proof}
We prove the claim by induction on $|R|$ and, for a fixed $j$, on
$|\Delta_1^{(j)}|$. Fix $j \in [1, 3]$ and assume that there is
either a unique source or a unique sink in $\Delta^{(l)}$ for each
$l \in [1, j - 1]$. Let $u$ and $v$ be the connecting vertices in
$\Delta$, and let $\alpha$ and $\beta$ be the arrows in
$\Delta^{(j)}$ adjacent to $u$ and $v$, respectively. The claim
follows by the arguments presented after
Lemma~\ref{lemm_shift_rel_group}, unless the following condition
(or its dual) is satisfied: $s \alpha = u$, $t \beta = v$, there
exists $\alpha' \in \Delta_1$ with $t \alpha' = u$ and $\alpha
\alpha' \in R$, and there exists $\beta' \in \Delta_1$ with $s
\beta' = v$ and $\beta' \beta \in R$. Assume the above condition
is satisfied. Put $x = s \beta$ and $v' = t \beta'$. If
$\Delta^{(j)}$ is not an equioriented linear quiver, then by
applying APR-coreflections, shifts of relations, and
Lemma~\ref{lemm_shift_last}, we may assume that there exists
$\gamma \in \Delta_1$ with $\gamma \neq \beta$ and $s \gamma = x$.
Put $y = t \gamma$.

Assume there exists $\gamma' \in \Delta_1$ with $s \gamma' = y$
and $\gamma' \gamma \in R$. Put $y' = t \gamma'$. If $(\Delta',
R')$ is the bound quiver obtained from $(\Delta, R)$ by applying
the APR-coreflections at $x$ and $y$, then $R' = R \setminus \{
\beta' \beta, \gamma' \gamma \} \cup \{ \gamma' \beta \}$ and
$\Delta'$ is obtained from $\Delta$ by replacing the subquiver
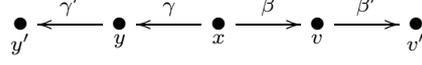
\[
\xymatrix{\vertexU{y'} & \vertexU{y} \ar[l]_{\gamma'} &
\vertexU{x} \ar[l]_\gamma \ar[r]^\beta & \vertexU{v}
\ar[r]^{\beta'} & \vertexU{v'}}
\]
by the quiver
\[
\vcenter
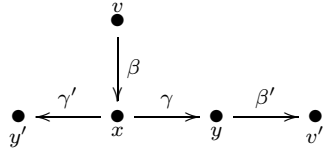
{\xymatrix{& \vertexD{v} \ar[d]^\beta \\ \vertexU{y'} &
\vertexU{x} \ar[l]_{\gamma'} \ar[r]^\gamma & \vertexU{y}
\ar[r]^{\beta'} & \vertexU{v'}}}.
\]
In particular, $|R'| < |R|$, hence the claim follows by induction
in this case.

Otherwise, if $(\Delta', R')$ is the bound quiver obtained from
$(\Delta, R)$ by applying the APR-coreflection at $x$, then $R' =
R \setminus \{ \beta' \beta \} \cup \{ \beta' \gamma \}$ and
$\Delta'$ is obtained from $\Delta$ by replacing the subquiver
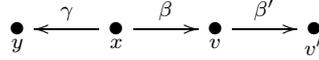
\[
\xymatrix{\vertexU{y} & \vertexU{x} \ar[l]_\gamma \ar[r]^\beta &
\vertexU{v} \ar[r]^{\beta'} & \vertexU{v'}}
\]
by the quiver
\[
\vcenter
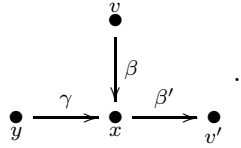
{\xymatrix{& \vertexD{v} \ar[d]^\beta \\
\vertexU{y} \ar[r]^\gamma & \vertexU{x} \ar[r]^{\beta'} &
\vertexU{v'}}}.
\]
Observe that if $l \in [1, j - 1]$ and there is no $\delta \in
\Delta_1^{(l)}$ with $s \delta = v$ and $\delta \neq \beta'$, then
there is either a unique source or a unique sink in
$\Delta'^{(l)}$. On the other hand, if there exists such an arrow,
then we may assume that there is either a unique source or a
unique sink in $\Delta'^{(l)}$, since $\beta$ is a free arrow in
$(\Delta', R')$. In particular, in both cases the claim follows
again by induction, since $|\Delta'^{(j)}_1| < |\Delta^{(j)}_1|$.
\end{proof}

\begin{step}
We may assume that if either $s \alpha = x = s \beta$ or $t \alpha
= x = t \beta$ for a connecting vertex $x$, $\alpha \in \Delta_{1,
\varepsilon_1} ^{(j_1)}$, and $\beta \in \Delta_{1,
\varepsilon_2}^{(j_2)}$, with $j_1 \neq j_2$ and $\varepsilon_1,
\varepsilon_2 \in \{ -, + \}$, then either
$R_{\varepsilon_1}^{(j_1)} = \varnothing$ or
$R_{\varepsilon_2}^{(j_2)} = \varnothing$.
\end{step}

\begin{proof}
Without loss of generality we may assume that $s \alpha = x = s
\beta$. If $R_{\varepsilon_1}^{(j_1)} \neq \varnothing$ or
$R_{\varepsilon_2}^{(j_2)} \neq \varnothing$, then by shifting
relations we may assume that there exist arrows $\alpha'$ and
$\beta'$ in $\Delta$ with $s \alpha' = t \alpha$, $s \beta' = t
\beta$, and $\alpha' \alpha, \beta' \beta \in R$. Let $\gamma$ be
the arrow in $\Delta$ with $t \gamma = x$. Without loss of
generality we may assume that $\alpha \gamma \in R$ and $\beta
\gamma \not \in R$. Put $y' = t \alpha'$, $z' = t \beta'$, and $u
= s \gamma$. If $(\Delta', R')$ is the bound quiver obtained from
$(\Delta, R)$ by applying the generalized APR-coreflections at $x$
and $y$, then $R' = R \setminus \{ \alpha' \alpha, \beta' \beta,
\alpha \gamma \} \cup \{ \alpha' \beta, \beta \gamma \}$ and
$\Delta'$ is obtained from $\Delta$ by replacing the subquiver
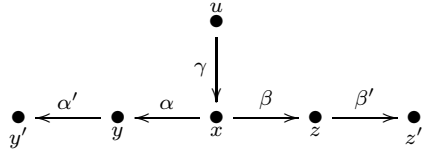
\[
\xymatrix{& & \vertexD{u} \ar[d]_\gamma \\
\vertexU{y'} & \vertexU{y} \ar[l]_{\alpha'} & \vertexU{x}
\ar[l]_\alpha \ar[r]^\beta & \vertexU{z} \ar[r]^{\beta'} &
\vertexU{z'}}
\]
by the quiver
\[
\vcenter
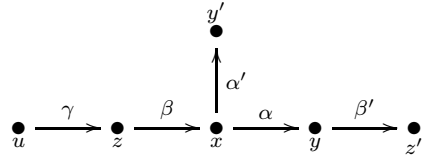
{\xymatrix{& & \vertexD{y'} \\
\vertexU{u} \ar[r]^\gamma & \vertexU{z} \ar[r]^\beta & \vertexU{x}
\ar[r]^\alpha \ar[u]_{\alpha'} & \vertexU{y} \ar[r]^{\beta'} &
\vertexU{z'}}}.
\]
In particular, $|R'^{(j_1)}_{\varepsilon_1}| <
|R^{(j_1)}_{\varepsilon_1}|$ and $|R'^{(j_2)}_{\varepsilon_2}| <
|R^{(j_2)}_{\varepsilon_2}|$, hence the claim follows by
induction.
\end{proof}

\begin{step}
We may assume that there exists $j \in [1, 3]$ such that
$\Delta^{(j)}$ is equioriented.
\end{step}

\begin{proof}
If the above condition is not satisfied, then without loss of
generality we may assume that
\begin{align*}
\Delta^{(1)} & = \xymatrix{\vertexU{u} & \vertexU{x_1}
\ar[l]_{\alpha_1} & & \vertexU{x_{p_1 - 1}} \ar[ll]|{\textstyle
\cdots} & \vertexU{x} \ar[l]_{\alpha_{p_1}} \ar[r]^{\alpha_{q_1}'}
& \vertexU{x_{q_1 - 1}'} \ar[rr]|{\textstyle \cdots} & &
\vertexU{x_1'} \ar[r]^{\alpha_1'} & \vertexU{v}},
\\
\Delta^{(2)} & = \xymatrix{\vertexU{u} \ar[r]^{\beta_{p_2}} &
\vertexU{y_{p_2 - 1}} \ar[rr]|{\textstyle \cdots} & &
\vertexU{y_1} \ar[r]^{\beta_1} & \vertexU{y} & \vertexU{y_1'}
\ar[l]_{\beta_1'} & & \vertexU{y_{q_2 - 1}'} \ar[ll]|{\textstyle
\cdots} & \vertexU{v} \ar[l]_{\beta_{q_2}'}},
\\
\intertext{and} %
\Delta^{(3)} & = \xymatrix{\vertexU{u} \ar[r]^{\gamma_{p_3}} &
\vertexU{z_{p_3 - 1}} \ar[rr]|{\textstyle \cdots} & &
\vertexU{z_1} \ar[r]^{\gamma_1} & \vertexU{z} & \vertexU{z_1'}
\ar[l]_{\gamma_1'} & & \vertexU{z_{q_3 - 1}'} \ar[ll]|{\textstyle
\cdots} & \vertexU{v} \ar[l]_{\gamma_{q_3}'}},
\end{align*}
for some $p_1, p_2, p_3, q_1, q_2, q_3 \in \bbN_+$. Moreover, we
may assume that $\beta_{p_2} \alpha_1 \in R$. Consequently, by
shifting relations we may assume that $\beta_i \beta_{i + 1} \in
R$ for all $i \in [1, p_2 - 1]$. There are two cases we have to
consider.

Assume first $\gamma_i \gamma_{i + 1} \not \in R$ for all $i \in
[1, p_3 - 1]$. If $(\Delta', R')$ is the bound quiver obtained
from $(\Delta, R)$ by applying the generalized APR-coreflection at
$u$ followed by the composition of the APR-coreflection at $y_i$
and the generalized APR-coreflection at $u$ for $i = p_2 - 1,
\ldots, 1$, then
\begin{align*}
R' & =
\begin{cases}
R \setminus \{ \beta_{p_2} \alpha_1, \beta_1 \beta_2 \} \cup \{
\gamma_{p_3} \alpha_1, \beta_{p_2} \gamma_{p_3} \} & p_2 > 1
\\
R \setminus \{ \beta_{p_2} \alpha_1 \} \cup \{ \gamma_{p_3}
\alpha_1 \} & p_2 = 1,
\end{cases}
\\
\Delta'^{(1)} & = \xymatrix{\vertexU{z_{p_3 - 1}} & \vertexU{x_1}
\ar[l]_{\alpha_1} & & \vertexU{x_{p_1 - 1}} \ar[ll]|{\textstyle
\cdots} & \vertexU{x} \ar[l]_{\alpha_{p_1}} \ar[r]^{\alpha_{q_1}'}
& \vertexU{x_{q_1 - 1}'} \ar[rr]|{\textstyle \cdots} & &
\vertexU{x_1'} \ar[r]^{\alpha_1'} & \vertexU{v}},
\\
\Delta'^{(2)} & = \makebox[\branch][l]{\xymatrix{\vertexU{z_{p_3 -
1}} \ar[r]^{\gamma_{p_3}} & \vertexU{y_{p_2 - 1}}
\ar[r]^{\beta_{p_2}} & \vertexU{y_{p_2 - 2}} \ar[rr]|{\textstyle
\cdots} & & \vertexU{u} & \vertexU{y} \ar[l]_{\beta_1} &
\vertexU{y_1'} \ar[l]_{\beta_1'} & & \vertexU{y_{q_2 - 1}'}
\ar[ll]|{\textstyle \cdots} & \vertexU{v} \ar[l]_{\beta_{q_2}'}},}
\\
\intertext{and} %
\Delta'^{(3)} & = \xymatrix{\vertexU{z_{p_3 - 1}}
\ar[rr]|{\textstyle \cdots} & & \vertexU{z_1} \ar[r]^{\gamma_1} &
\vertexU{z} & \vertexU{z_1'} \ar[l]_{\gamma_1'} & &
\vertexU{z_{q_3 - 1}'} \ar[ll]|{\textstyle \cdots} & \vertexU{v}
\ar[l]_{\gamma_{q_3}'}},
\end{align*}
where $z_{p_3 - 1} = z$ if $p_3 = 1$. Consequently, the claim
follows by an easy induction.

Assume now that there exists $i \in [1, p_3 - 1]$ such that
$\gamma_i \gamma_{i + 1} \in R$. Consequently $p_2 = 1$. Moreover,
by shifting relations we may assume that $\gamma_{p_3 - 1}
\gamma_{p_3} \in R$. If $(\Delta', R')$ is the bound quiver
obtained from $(\Delta, R)$ by applying the generalized
APR-coreflection at $u$, then $R' = R \setminus \{ \beta_1
\alpha_1, \gamma_{p_3 - 1} \gamma_{p_3} \} \cup \{ \gamma_{p_3}
\alpha_1, \gamma_{p_3 - 1} \beta_1 \}$,
\begin{align*}
\Delta'^{(1)} & = \makebox[\branch][l]{\xymatrix{\vertexU{u} &
\vertexU{z_{p_3 - 1}} \ar[l]_{\gamma_{p_3}} & \vertexU{x_1}
\ar[l]_{\alpha_1} & & \vertexU{x_{p_1 - 1}} \ar[ll]|{\textstyle
\cdots} & \vertexU{x} \ar[l]_{\alpha_{p_1}} \ar[r]^{\alpha_{q_1}'}
& \vertexU{x_{q_1 - 1}'} \ar[rr]|{\textstyle \cdots} & &
\vertexU{x_1'} \ar[r]^{\alpha_1'} & \vertexU{v}},}
\\
\Delta'^{(2)} & = \xymatrix{ \vertexU{u} & \vertexU{y}
\ar[l]_{\beta_1} & \vertexU{y_1'} \ar[l]_{\beta_1'} & &
\vertexU{y_{q_2 - 1}'} \ar[ll]|{\textstyle \cdots} & \vertexU{u}
\ar[l]_{\beta_{q_2}'}},
\\
\intertext{and} %
\Delta'^{(3)} & = \xymatrix{\vertexU{u} \ar[r]^{\gamma_{p_3 - 1}}
& \vertexU{z_{p_3 - 2}} \ar[rr]|{\textstyle \cdots} & &
\vertexU{z} & \vertexU{z_1'} \ar[l]_{\gamma_1'} & &
\vertexU{z_{q_3 - 1}'} \ar[ll]|{\textstyle \cdots} & \vertexU{u}
\ar[l]_{\gamma_{q_3}'}},
\end{align*}
thus the claim follows.
\end{proof}

\begin{step}
We may assume that there is at most one $j \in [1, 3]$ such that
$\Delta^{(j)}$ is not equioriented.
\end{step}

\begin{proof}
If the above condition is not satisfied, then without loss of
generality we may assume that
\begin{align*}
\Delta^{(1)} & = \xymatrix{\vertexU{u} & \vertexU{x_1}
\ar[l]_{\alpha_1} & & \vertexU{x_{p_1 - 1}} \ar[ll]|{\textstyle
\cdots} & \vertexU{x} \ar[l]_{\alpha_{p_1}} \ar[r]^{\alpha_{q_1}'}
& \vertexU{x_{q_1 - 1}'} \ar[rr]|{\textstyle \cdots} & &
\vertexU{x_1'} \ar[r]^{\alpha_1'} & \vertexU{v}},
\\
\Delta^{(2)} & = \xymatrix{\vertexU{u} \ar[r]^{\beta_{p_2}} &
\vertexU{y_{p_2 - 1}} \ar[rr]|{\textstyle \cdots} & &
\vertexU{y_1} \ar[r]^{\beta_1} & \vertexU{y} & \vertexU{y_1'}
\ar[l]_{\beta_1'} & & \vertexU{y_{q_2 - 1}'} \ar[ll]|{\textstyle
\cdots} & \vertexU{v} \ar[l]_{\beta_{q_2}'}},
\\
\intertext{and} %
\Delta^{(3)} & = \xymatrix{\vertexU{u} \ar[r]^{\gamma_{p_3}} &
\vertexU{z_{p_3 - 1}} \ar[rr]|{\textstyle \cdots} & &
\vertexU{z_1} \ar[r]^{\gamma_1} & \vertexU{v}},
\end{align*}
for some $p_1, p_2, p_3, q_1, q_2 \in \bbN_+$. In this proof we
will again number the cases. Up to symmetry, there are three main
cases we have to consider: either $\beta_{p_2} \alpha_1 \in R$ and
$\beta_{q_2}' \alpha_1' \in R$, or $\beta_{p_2} \alpha_1 \in R$
and $\beta_{q_2}' \gamma_1 \in R$, or $\gamma_{p_3} \alpha_1 \in
R$ and $\beta_{q_2}' \gamma_1 \in R$.

(1)~Assume $\beta_{p_2} \alpha_1 \in R$ and $\beta_{q_2}'
\alpha_1' \in R$. In this case we may apply the same arguments as
in the proof of the previous step. Note however that if $\gamma_i
\gamma_{i + 1} \not \in R$ for all $i \in [1, p_3 - 1]$, then we
obtain a gentle bound quiver whose bound quiver algebra is
tilting-cotilting equivalent to $\Lambda_0' (p, r)$ for some $p
\in \bbN_+$ and $r \in [0, p - 1]$ according to
Proposition~\ref{prop_Lambda0}.

(2)~Assume that $\beta_{p_2} \alpha_1 \in R$ and $\beta_{q_2}'
\gamma_1 \in R$.

(2.1)~Assume that $\alpha_i' \alpha_{i + 1}' \not \in R$ for all
$i \in [1, q_1 - 1]$ and $\beta_i' \beta_{i + 1}' \not \in R$ for
all $i \in [1, q_2 - 1]$. By shifting the relation $\beta_{q_2}'
\gamma_1$ to the left we may assume that $q_2 = 1$.

(2.1.1)~Assume that $\beta_i \beta_{i + 1} \not \in R$ for all $i
\in [1, p_2 - 1]$. By shifting the relation $\beta_{p_2} \alpha_1$
to the left we may assume that $p_2 = 1$. Consequently, the path
algebra of the bound quiver obtained from $(\Delta, R)$ by
application of the APR-reflections at $y$, $v$, $x_1'$, \ldots,
$x_{q_1 - 1}'$ is easily seen to be tilting-cotilting equivalent
to $\Lambda_2 (p_3 + 1, p_1 + 1, q_1, r_3, r_1)$, where $r_1$ is
the number of $i \in [1, p_1 - 1]$ such that $\alpha_i \alpha_{i +
1} \in R$ and $r_3$ is the number of $i \in [1, p_3 - 1]$ such
that $\gamma_i \gamma_{i + 1} \in R$.

(2.1.2)~Assume that there exists $i \in [1, p_2 - 1]$ such that
$\beta_i \beta_{i + 1} \in R$. By shifting relations we may assume
that $\beta_1 \beta_2 \in R$. If $(\Delta', R')$ is the bound
quiver obtained from $(\Delta, R)$ by applying the APR-reflections
at $y$, $y_1$, $v$, $x_1'$, \ldots, $x_{q_1 - 1}'$, then $R' = R
\setminus \{ \beta_1 \beta_2, \beta_1' \gamma_1 \} \cup \{
\alpha_{q_1}' \beta_2 \}$,
\begin{align*}
\Delta'^{(1)} & = \xymatrix{\vertexU{u} & \vertexU{x_1}
\ar[l]_{\alpha_1} & & \vertexU{x_{p_1 - 1}} \ar[ll]|{\textstyle
\cdots} & \vertexU{x} \ar[l]_{\alpha_{p_1}} & \vertexU{x_{q_1 -
1}'} \ar[l]_{\alpha_{q_1}'}},
\\
\Delta'^{(2)} & = \xymatrix{\vertexU{u} \ar[rr]|{\textstyle
\cdots} & & \vertexU{y_2} \ar[r]^{\beta_2} & \vertexU{x_{q_1 -
1}'}},
\\
\intertext{and} %
\Delta'^{(3)} & = \makebox[\branch][l]{\xymatrix{\vertexU{u}
\ar[r]^{\gamma_{p_3}} & \vertexU{z_{p_3 - 1}} \ar[rr]|{\textstyle
\cdots} & & \vertexU{z_1} \ar[r]^{\gamma_1} & \vertexU{y_1}
\ar[r]^{\beta_1} & \vertexU{y} & \vertexU{v} \ar[l]_{\beta_1'} &
\vertexU{x_1'} \ar[l]_{\alpha_1'} & & \ar[ll]|{\textstyle \cdots}
\vertexU{x_{q_1 - 1}'}},}
\\
\end{align*}
where $x_{q_1 - 1}' = v$ if $q_1 = 1$, hence the claim follows.

(2.2)~Assume that there exists $i \in [1, q_2 - 1]$ such that
$\beta_i' \beta_{i + 1}' \in R$. By shifting relations we may
assume that $\beta_1' \beta_2' \in R$. Moreover, this condition
implies that $\beta_i \beta_{i + 1} \not \in R$ for all $i \in [1,
p_2 - 1]$. By shifting the relation $\beta_{p_2} \alpha_1$ to the
left we may assume $p_2 = 1$. If $(\Delta', R')$ is the bound
quiver obtained from $(\Delta, R)$ by applying the APR-reflections
at $y$ and $y_1'$, then $R' = R \setminus \{ \beta_1 \alpha_1,
\beta_1' \beta_2' \} \cup \{ \beta_1 \beta_2' \}$,
\begin{align*}
\Delta'^{(1)} & = \makebox[\branch][l]{\xymatrix{\vertexU{y} &
\vertexU{y_1'} \ar[l]_{\beta_1'} & \vertexU{x_1} \ar[l]_{\alpha_1}
& & \vertexU{x_{p_1 - 1}} \ar[ll]|{\textstyle \cdots} &
\vertexU{x} \ar[l]_{\alpha_{p_1}} \ar[r]^{\alpha_{q_1}'} &
\vertexU{x_{q_1 - 1}'} \ar[rr]|{\textstyle \cdots} & &
\vertexU{x_1'} \ar[r]^{\alpha_1'} & \vertexU{v}},}
\\
\Delta'^{(2)} & = \xymatrix{\vertexU{y} & \vertexU{y_2'}
\ar[l]_{\beta_2'} & & \vertexU{v} \ar[ll]|{\textstyle \cdots}},
\\
\intertext{and} %
\Delta'^{(3)} & = \xymatrix{\vertexU{y} \ar[r]^{\beta_1} &
\vertexU{u} \ar[r]^{\gamma_{p_3}} & \vertexU{z_{p_3 - 1}}
\ar[rr]|{\textstyle \cdots} & & \vertexU{z_1} \ar[r]^{\gamma_1} &
\vertexU{v}},
\end{align*}
hence the claim follows again.

(2.3)~Assume that $\beta_i' \beta_{i + 1}' \not \in R$ for all $i
\in [1, q_2 - 1]$ and there exists $i \in [1, q_1 - 1]$ such that
$\alpha_i' \alpha_{i + 1}' \in R$. Observe that in this case
$\alpha_i \alpha_{i + 1} \not \in R$ for all $i \in [1, p_1 - 1]$,
hence by shifting the relation $\beta_{p_2} \alpha_1$ to the right
we may assume that $p_1 = 1$. Similarly, $\gamma_i \gamma_{i + 1}
\not \in R$ for all $i \in [1, p_3 - 1]$.

(2.3.1)~Assume that there exists $i \in [1, p_2 - 1]$ such that
$\beta_i \beta_{i + 1} \in R$, then by shifting relations we may
assume that $\beta_1 \beta_2 \in R$. Moreover, by shifting the
relation $\beta_{q_2}' \gamma_1$ to the left we may assume that
$q_2 = 1$. Additionally, by shifting relations we may assume that
$\alpha_1' \alpha_2' \in R$. If $(\Delta', R')$ is the bound
quiver obtained from $(\Delta, R)$ by applying the HW-reflection
at $y$ followed by the APR-coreflection at $y$ and the
APR-reflections at $v$, $z_1$, \ldots, $z_{p_3 - 1}$, then $R' = R
\setminus \{ \alpha_1' \alpha_2', \beta_1 \beta_2, \beta_1'
\gamma_1 \} \cup \{ \gamma_{p_3} \alpha_2' \}$,
\begin{align*}
\Delta'^{(1)} & = \xymatrix{\vertexU{u} & \vertexU{x}
\ar[l]_{\alpha_1} \ar[rr]|{\textstyle \cdots} & & \vertexU{x_2'}
\ar[r]^{\alpha_2'} & \vertexU{z_{p_3 - 1}}},
\\
\Delta'^{(2)} & = \makebox[\branch][l]{\xymatrix{\vertexU{u}
\ar[r]^{\beta_{p_2}} & \vertexU{y_{p_2 - 1}} \ar[rr]|{\textstyle
\cdots} && \vertexU{y_1} \ar[r]^{\beta_1} & \vertexU{y} &
\vertexU{x_1'} \ar[l]_{\beta_1'} & \vertexU{v} \ar[l]_{\alpha_1'}
& \vertexU{z_1} \ar[l]_{\gamma_1} & & \vertexU{z_{p_3 - 1}}
\ar[ll]|{\textstyle \cdots}},}
\\
\intertext{and} %
\Delta'^{(3)} & = \xymatrix{\vertexU{u} & \vertexU{z_{p_3 - 1}}
\ar[l]_{\gamma_{p_3}}},
\end{align*}
where $z_{p_3 - 1} = v$ if $p_3 = 1$. In particular, we reduce the
proof to the situation dual either to~(2.1) or to~(2.2).

(2.3.2)~Assume that $\beta_i \beta_{i + 1} \not \in R$ for all $i
\in [1, p_2 - 1]$. By shifting the relation $\beta_{q_2}'
\gamma_1$ to the right we may assume that $p_3 = 1$. Additionally,
by shifting relations we may assume that $\alpha_{q_1 - 1}'
\alpha_{q_1}' \in R$. The bound quiver algebra of the bound quiver
obtained from $(\Delta, R)$ by applying the APR-coreflections at
$x$, $x_{q_1 - 1}'$, $u$, $x$, $x_{q_1 - 1}'$, $y_{p_2 - 1}$,
\ldots, $y_1$ is easily seen to be titling-cotilting equivalent to
$\Lambda_2 (q_2 + 1, q_1, p_2 + 1, 0, r_1' - 1)$, where $r_1'$ is
the number of $i \in [1, q_1 - 1]$ such that $\alpha_i' \alpha_{i
+ 1}' \in R$.

(3)~Assume that $\gamma_{p_3} \alpha_1 \in R$ and $\beta_{q_2}'
\gamma_1 \in R$.

(3.1)~Assume that there exists $i \in [1, p_2 - 1]$ such that
$\beta_i \beta_{i + 1}  \in R$. By shifting relations we may
assume that $\beta_{p_2 - 1} \beta_{p_2} \in R$. Since in this
case $\gamma_i \gamma_{i + 1} \not \in R$ for all $i \in [1, p_3 -
1]$, we may assume, by shifting the relation $\gamma_{p_3}
\alpha_1$ to the left, that $p_3 = 1$. Consequently, the bound
quiver algebra of the bound quiver obtained from $(\Delta, R)$ by
applying the generalized APR-coreflection at $u$ is
tilting-cotilting equivalent to $\Lambda_0' (p, r)$ for some $p
\in \bbN_+$ and $r \in [0, p - 1]$ according to
Proposition~\ref{prop_Lambda0}.

(3.2)~Assume that $\beta_i \beta_{i + 1} \not \in R$ for all $i
\in [1, p_2 - 1]$. By shifting relations we may also assume that
$\gamma_i \gamma_{i + 1} \in R$ for all $i \in [1, p_3 - 1]$. If
$(\Delta', R')$ is the bound quiver obtained from $(\Delta, R)$ by
applying the APR-coreflection at $u$ followed by the composition
of the APR-coreflection at $z_i$ and the generalized
APR-coreflection at $u$ for $i = p_3 - 1, \ldots, 1$, then
\begin{align*}
R' & =
\begin{cases}
R \setminus \{ \gamma_{p_3} \alpha_1, \gamma_1 \gamma_2,
\beta_{q_2}' \gamma_1 \} \cup \{ \beta_{p_2} \alpha_1,
\gamma_{p_3} \beta_{p_2}, \beta_{q_2}' \gamma_2 \} & p_3 > 1,
\\
R \setminus \{ \gamma_1 \alpha_1, \beta_{q_2}' \gamma_1 \} \cup \{
\beta_{p_2} \alpha_1, \beta_{q_2}' \beta_{p_2} \} & p_3 = 1,
\end{cases}
\\
\Delta'^{(1)} & = \makebox[\branch][l]{\xymatrix{\vertexU{y_{p_2 -
1}} & \vertexU{x_1} \ar[l]_{\alpha_1} & & \vertexU{x_{p_1 - 1}}
\ar[ll]|{\textstyle \cdots} & \vertexU{x} \ar[l]_{\alpha_{p_1}}
\ar[r]^{\alpha_{q_1}'} & \vertexU{x_{q_1 - 1}'}
\ar[rr]|{\textstyle \cdots} & & \vertexU{x_1'} \ar[r]^{\alpha_1'}
& \vertexU{v} \ar[r]^{\gamma_{p_3}} & \vertexU{u}},}
\\
\Delta'^{(2)} & = \xymatrix{\vertexU{y_{p_2 - 1}}
\ar[rr]|{\textstyle \cdots} & & \vertexU{y_1} \ar[r]^{\beta_1} &
\vertexU{y} & \vertexU{y_1'} \ar[l]_{\beta_1'} & & \vertexU{y_{q_2
-  1}'} \ar[ll]|{\textstyle \cdots} & \vertexU{u}
\ar[l]_{\beta_{q_2}'}}
\\
\intertext{and} %
\Delta'^{(3)} & = \xymatrix{\vertexU{y_{p_2 - 1}}
\ar[r]^{\beta_{p_2}} & \vertexU{z_{p_3 - 1}}\ar[rr]|{\textstyle
\cdots} & &  \vertexU{z_1} \ar[r]^{\gamma_2} & \vertexU{u}},
\end{align*}
where $y_{p_2 - 1} = y$ if $p_2 = 1$. Consequently, the claim
follows by induction.
\end{proof}

For $p_1, p_2, p_3 \in \bbN_+$, $p_2 \geq 2$, $r_1 \in [0, p_1 -
1]$, and $r_2 \in [1, p_2 - 1]$, let $\Lambda_2' (p_1, p_2, p_3,
r_1, r_2)$ be the algebra of the quiver
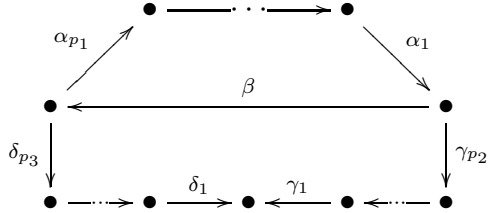
\[
\xymatrix{& \bullet \ar[rr]|{\textstyle \cdots} & & \bullet
\ar[rd]^{\alpha_1} \\ \bullet \ar[ru]^{\alpha_{p_1}}
\ar[d]_{\delta_{p_3}} & & & & \bullet \ar[d]^{\gamma_{p_2}}
\ar[llll]_\beta \\ \bullet \ar[r]|\cdots & \bullet
\ar[r]^{\delta_1} & \bullet & \bullet \ar[l]_{\gamma_1} & \bullet
\ar[l]|\cdots}
\]
bound by $\alpha_i \alpha_{i + 1}$ for $i \in [p_1 - r_1, p_1 -
1]$, $\alpha_{p_1} \beta$, $\beta \alpha_1$, $\gamma_i \gamma_{i +
1}$ for $i \in [1, r_2]$. Observe that $\Lambda_2' (p_1, p_2, p_3,
r_1, r_2)$ is tilting-cotilting equivalent to $\Lambda_2 (p_2, p_1
+ 1, p_3, r_2 - 1, r_1 + 1)$. Indeed, it is enough to apply the
HW-reflection at $x_i$ followed by the APR-coreflection at $x_i$
for $i = 1, \ldots, p_3$, where $x_i = t \delta_i$ for $i \in [1,
p_3]$.

\begin{step}
We may assume that $\Delta^{(j)}$ is equioriented for each $j \in
[1, 3]$.
\end{step}

\begin{proof}
Suppose that there exists $j \in [1, 3]$ such that $\Delta^{(j)}$
is not equioriented. Without loss of generality we may assume that
\begin{align*}
\Delta^{(1)} & = \xymatrix{\vertexU{u} & \vertexU{x_1}
\ar[l]_{\alpha_1} & & \vertexU{x_{p_1 - 1}} \ar[ll]|{\textstyle
\cdots} & \vertexU{x} \ar[l]_{\alpha_{p_1}} \ar[r]^{\alpha_{q_1}'}
& \vertexU{x_{q_1 - 1}'} \ar[rr]|{\textstyle \cdots} & &
\vertexU{x_1'} \ar[r]^{\alpha_1'} & \vertexU{v}},
\\
\Delta^{(2)} & = \xymatrix{\vertexU{u} \ar[r]^{\beta_{p_2}} &
\vertexU{y_{p_2 - 1}} \ar[rr]|{\textstyle \cdots} & &
\vertexU{y_1} \ar[r]^{\beta_1} & \vertexU{v}},
\\
\intertext{and} %
\Delta^{(3)} & = \xymatrix{\vertexU{u} & \vertexU{z_1}
\ar[l]_{\gamma_1} & & \vertexU{z_{p_3 - 1}} \ar[ll]|{\textstyle
\cdots} & \vertexU{v} \ar[l]_{\gamma_{p_3}}},
\end{align*}
for some $p_1, p_2, p_3, q_1 \in \bbN_+$. We may additionally
assume that $\alpha_{i + 1} \alpha_i \not \in R$ for all $i \in
[1, p_1 - 1]$. Let $r_1'$ be the number of $i \in [1, q_1 - 1]$
such that $\alpha_i' \alpha_{i + 1}' \in R$, let $r_2$ be the
number of $i \in [1, p_2 - 1]$ such that $\beta_i \beta_{i + 1}
\in R$, and let $r_3$ be the number of $i \in [1, p_3 - 1]$ such
that $\gamma_i \gamma_{i + 1} \in R$. Observe that by symmetry we
may assume that $r_1' > 0$ if $\gamma_{p_3} \alpha_1' \in R$ and
$\beta_{p_2} \alpha_1 \not \in R$.

Assume first that $\beta_{p_2} \alpha_1 \in R$. In this case by
shifting the relation $\beta_{p_2} \alpha_1$ to the right we may
assume that $p_1 = 1$. Observe that either $r_1' = 0$ or $r_2 =
0$.

If $r_3 \geq 1$, then by shifting relations we may assume that
$\gamma_1 \gamma_2 \in R$. If $(\Delta', R')$ is the bound quiver
obtained from $(\Delta, R)$ by applying the generalized
APR-reflection at $u$, then $R' = R \setminus \{ \beta_{p_2}
\alpha_1, \gamma_1 \gamma_2 \} \cup \{ \alpha_1 \gamma_2,
\beta_{p_2} \gamma_1 \}$,
\begin{align*}
\Delta'^{(1)} & = \xymatrix{\vertexU{u} \ar[r]^{\alpha_1} &
\vertexU{x} \ar[r]^{\alpha_{q_1}'} & \vertexU{x_{q_1 - 1}'}
\ar[rr]|{\textstyle \cdots} & & \vertexU{x_1'} \ar[r]^{\alpha_1'}
& \vertexU{v}},
\\
\Delta'^{(2)} & = \xymatrix{\vertexU{u} \ar[r]^{\gamma_1} &
\vertexU{z_1} \ar[r]^{\beta_{p_2}} & \vertexU{y_{p_2 - 1}}
\ar[rr]|{\textstyle \cdots} & & \vertexU{y_1} \ar[r]^{\beta_1} &
\vertexU{v}},
\\
\intertext{and} %
\Delta'^{(3)} & = \xymatrix{\vertexU{u} & \vertexU{z_2}
\ar[l]_{\gamma_2} & & \vertexU{v} \ar[ll]|{\textstyle \cdots}},
\end{align*}
hence the claim follows in this case.

Assume now that $r_3 = 0$. There are two additional possibilities
in this case. If $\gamma_{p_3} \alpha_1' \in R$, then $r_2 \geq 1$
(since $(\Delta, R)$ is a bound quiver). Consequently, $r_1' = 0$
and we have the situation symmetric to the previous one. On the
other hand, if $\gamma_{p_3} \beta_1 \in R$, then by shifting the
relation $\gamma_{p_3} \beta_1$ to the left we may assume that
$p_3 = 1$. Consequently, if $(\Delta', R')$ is the bound quiver
obtained from $(\Delta, R)$ by applying the generalized
APR-reflection at $u$, then $R' = R \setminus \{ \beta_{p_2}
\alpha_1, \gamma_1 \beta_1 \} \cup \{ \alpha_1 \beta_1,
\beta_{p_2} \gamma_1 \}$,
\begin{align*}
\Delta'^{(1)} & = \xymatrix{\vertexU{u} \ar[r]^{\alpha_1} &
\vertexU{x} \ar[r]^{\alpha_{q_1}'} & \vertexU{x_{q_1 - 1}'}
\ar[rr]|{\textstyle \cdots} & & \vertexU{x_1'} \ar[r]^{\alpha_1'}
& \vertexU{v}},
\\
\Delta'^{(2)} & = \xymatrix{\vertexU{u} & \vertexU{y_1}
\ar[l]_{\beta_1} & & \vertexU{y_{p_2 - 1}} \ar[ll]|{\textstyle
\cdots} & \vertexU{v} \ar[l]_{\beta_{p_2}}},
\\
\intertext{and} %
\Delta'^{(3)} & = \xymatrix{\vertexU{u} \ar[r]^{\gamma_1} &
\vertexU{v}},
\end{align*}
hence the claim follows.

Assume now that $\beta_{p_2} \gamma_1 \in R$. If $\gamma_{p_3}
\beta_1 \in R$, then it is follows easily that the bound quiver
algebra of $(\Delta, R)$ is tilting-cotilting equivalent either to
$\Lambda_2 (p_2 + p_3 - r_2 - 1, r_2 + 1, q_1, p_1, r_3)^{\op}$ if
$r_1' = 0$, or to $\Lambda_2' (p_2 + p_3 - 1, q_1, p_1, r_3,
r_1')^{\op}$ if $r_1' \geq 1$. Since $\Lambda_2' (p_2 + p_3 - 1,
q_1, p_1, r_3, r_1')$ is tilting-cotilting equivalent to
$\Lambda_2 (q_1, p_2 + p_3, p_1, r_1' - 1, r_3 + 1)$, thus we may
assume that $\gamma_{p_3} \alpha_1' \in R$. Consequently, by
shifting relations we may assume that $\alpha_i' \alpha_{i + 1}'
\in R$ for all $i \in [1, q_1 - 1]$. Recall that $q_1 > 1$ in this
case. If $(\Delta', R')$ is the bound quiver obtained from
$(\Delta, R)$ by applying the APR-coreflection at $x$ followed by
the composition of the HW-coreflection at $x_{q_1 - 1}'$ and the
APR-reflection at $x_{q_1 - 1}'$ applied $q_1 - 1$ times, then $R'
= R \setminus \{ \alpha_{q_1 - 1}' \alpha_{q_1}' \} \cup \{
\alpha_{q_1 - 1}' \alpha_{p_1} \}$,
\begin{align*}
\Delta'^{(1)} & =
\begin{cases}
\xymatrix{\vertexU{u} & \vertexU{x_1} \ar[l]_{\alpha_1} & &
\vertexU{x_{p_1 - 1}} \ar[ll]|{\textstyle \cdots}
\ar[r]^{\alpha_{p_1}} & \vertexU{x} \ar[r]^{\alpha_{q_1 - 1}'} &
\vertexU{x_{q_1 - 2}'} \ar[rr]|{\textstyle \cdots} & &
\vertexU{v}} & p_1 > 1,
\\
\xymatrix{\vertexU{u} \ar[r]^{\alpha_1} & \vertexU{x}
\ar[r]^{\alpha_{q_1 - 1}'} & \vertexU{x_{q_1 - 2}'}
\ar[rr]|{\textstyle \cdots} & & \vertexU{v}} & p_1 = 1,
\end{cases}
\\
\Delta'^{(2)} & = \xymatrix{\vertexU{u} \ar[r]^{\beta_1} &
\vertexU{y_1} \ar[rr]|{\textstyle \cdots} & & \vertexU{y_{p_2 -
1}} \ar[r]^{\beta_{p_2}} & \vertexU{x_{q_1 - 1}'}
\ar[r]^{\alpha_{q_1}'} & \vertexU{v}},
\\
\intertext{and} %
\Delta'^{(3)} & = \xymatrix{\vertexU{u} & \vertexU{z_1}
\ar[l]_{\gamma_1} & & \vertexU{z_{p_3 - 1}} \ar[ll]|{\textstyle
\cdots} & \vertexU{v} \ar[l]_{\gamma_{p_3}}},
\end{align*}
thus the claim follows by induction.
\end{proof}

We may prove now Proposition~\ref{prop_mainprop}. According to our
considerations we may assume that $(\Delta, R)$ is proper,
\begin{align*}
\Delta^{(1)} & = \xymatrix{\vertexU{u} & \vertexU{x_1}
\ar[l]_{\alpha_1} & & \vertexU{x_{p_1 - 1}} \ar[ll]|{\textstyle
\cdots} & \vertexU{v} \ar[l]_{\alpha_{p_1}}},
\\
\Delta^{(2)} & = \xymatrix{\vertexU{u} \ar[r]^{\beta_{p_2}} &
\vertexU{y_{p_2 - 1}} \ar[rr]|{\textstyle \cdots} & &
\vertexU{y_1} \ar[r]^{\beta_1} & \vertexU{v}},
\\
\intertext{and} %
\Delta^{(3)} & = \xymatrix{\vertexU{u} & \vertexU{z_1}
\ar[l]_{\gamma_1} & & \vertexU{z_{p_3 - 1}} \ar[ll]|{\textstyle
\cdots} & \vertexU{v} \ar[l]_{\gamma_{p_3}}},
\end{align*}
for some $p_1, p_2, p_3 \in \bbN_+$. Moreover, we may
additionally assume that $\alpha_i \alpha_{i + 1} \not \in R$ for
all $i \in [1, p_1 - 1]$. Let $r_2$ be the number of $i \in [1,
p_2 - 1]$ such that $\beta_i \beta_{i + 1} \in R$ and let $r_3$ be
the number of $i \in [1, p_3 - 1]$ such that $\gamma_i \gamma_{i +
1} \in R$.

Observe that if either $\beta_{p_2} \alpha_1 \in R$ and
$\gamma_{p_3} \beta_1 \in R$ or $\beta_{p_2} \gamma_1 \in R$ and
$\alpha_{p_1} \beta_1 \in R$, then the claim follows from
Proposition~\ref{prop_Lambda0prim}, thus we have to consider two
remaining cases.

Assume first that $\beta_{p_2} \alpha_1 \in R$ and $\alpha_{p_1}
\beta_1 \in R$. In this case by shifting the relation $\beta_{p_2}
\alpha_1$ to the right we may assume that $p_1 = 1$. If $r_3 = 0$,
then the bound quiver algebra of $(\Delta, R)$ is
tilting-cotilting equivalent to the $\Lambda_1 (p_2, 1, p_3, 0,
r_2)$ (observe that $r_2 \geq 1$ since $(\Delta, R)$ is a bound
quiver). On the other hand, if there exists $i \in [1, p_3 - 1]$
such that $\gamma_i \gamma_{i + 1} \in R$, then by shifting
relations we may assume that $\gamma_1 \gamma_2 \in R$.
Consequently, the bound quiver algebra of the bound quiver
obtained from $(\Delta, R)$ by applying the generalized
APR-reflection at $u$ is tilting-cotilting equivalent to
$\Lambda_2 (p_2 + 1, p_3, 0, r_2 + 1, r_3 - 1)$ and this finishes
the proof in this case.

Assume now that $\beta_{p_2} \gamma_1 \in R$ and $\gamma_{p_3}
\beta_1 \in R$. In this case it follows by shifting relations that
the bound quiver algebra of $(\Delta, R)$ is tilting-cotilting
equivalent to $\Lambda_1 (p_2 + p_3 - r_3 - 1, r_3 + 1, p_1, 0,
r_2)$ (again $r_2 \geq 1$ since $(\Delta, R)$ is a bound quiver)
and this finishes the proof.

\section{Minimality of the list} \label{sect_min}

In this section we prove that different algebras from the list in
Theorem~\ref{theo_nondeg} are not derived equivalent. We also
check that the algebras listed in Theorem~\ref{theo_nondeg} are
nondegenerate, while the algebras listed in Theorem~\ref{theo_deg}
are degenerate. A tool used in order to distinguish between
derived equivalence classes of these algebras will be the derived
invariant introduced by Avella-Alaminos and Geiss
in~\cite{AAGei2006}.

Let $(\Delta, R)$ be a gentle quiver. By a permitted thread in
$(\Delta, R)$ we mean either a maximal path in $(\Delta, R)$ or $x
\in \Delta_0$ such that there is at most one arrow $\alpha$ with
$s \alpha = x$, there is at most one arrow $\beta$ with $t \beta =
x$, and $\alpha \beta \not \in R$ for all $\alpha, \beta \in
\Delta_1$ with $s \alpha = x = t \beta$. Similarly, we define
notion of a forbidden thread in $(\Delta, R)$. Namely, first we
say that by an anti-path in $(\Delta, R)$ we mean a path $\alpha_1
\cdots \alpha_n$ in $\Delta$ such that $\alpha_i \alpha_{i + 1}
\in R$ for all $i \in [1, n - 1]$. In particular, every trivial
path is an anti-path. By a forbidden thread we mean either a
maximal anti-path in $(\Delta, R)$ or $x \in \Delta_0$ such that
there is at most one arrow $\alpha$ with $s \alpha = x$, there is
at most one arrow $\beta$ with $t \beta = x$, and $\alpha \beta
\in R$ for all $\alpha, \beta \in \Delta_1$ with $s \alpha = x = t
\beta$.

By a characteristic sequence in a gentle bound quiver $(\Delta,
R)$ we mean a sequence $(\sigma_i, \tau_i)_{i \in \bbZ}$ of
permitted threads $\sigma_i$, $i \in \bbZ$, and forbidden threads
$\tau_i$, $i \in \bbZ$, such that for each $i \in \bbZ$ the
following conditions are satisfied:
\begin{enumerate}

\item
$t \tau_i = t \sigma_i$ and $s \sigma_{i + 1} = s \tau_i$,

\item
if $\sigma_i = x = \tau_i$ for $x \in \Delta_0$ then $\sigma_{i +
1} \neq x$, unless $\Delta_1 = \varnothing$,

\item
if $\tau_i = x = \sigma_{i + 1}$ for $x \in \Delta_0$, then
$\tau_{i + 1} \neq x$, unless $\Delta_1 = \varnothing$,

\item
if neither $\sigma_i$ nor $\tau_i$ is a trivial path, then the
terminating arrow of $\tau_i$ differs from the terminating arrow
of $\sigma_i$,

\item
if neither $\tau_i$ nor $\sigma_{i + 1}$ is a trivial path, then
the starting arrow of $\sigma_{i + 1}$ differs from the starting
arrow of $\tau_i$.

\end{enumerate}
We identify characteristic sequences $(\sigma_i, \tau_i)_{i \in
\bbZ}$ and $(\sigma_i', \tau_i')_{i \in \bbZ}$ if there exists $l
\in \bbZ$ such that $\sigma_i = \sigma_{i + l}'$ and $\tau_i =
\tau_{i + l}'$ for all $i \in \bbZ$. By the type of the
characteristic sequence $(\sigma_i, \tau_i)_{i \in \bbZ}$ we mean
a pair $(n, m) \in \bbN \times \bbN$ defined by $n = \min \{ l \in
\bbN_+ \mid \sigma_l = \sigma_0 \}$ and $m = \sum_{i \in [1, n]}
\ell (\tau_i)$. In the above situation we also write $(\sigma_1,
\tau_1, \cdots, \sigma_n, \tau_n)$ instead of $(\sigma_i,
\tau_i)_{i \in \bbZ}$. Additionally, we also call a sequence
$(\alpha_i)_{i \in \bbZ}$ of arrows in $\Delta$ a characteristic
sequence if $s \alpha_i = t \alpha_{i + 1}$ and $\alpha_i
\alpha_{i + 1} \in R$ for all $i \in \bbZ$. Again we identify
sequences $(\alpha_i)_{i \in \bbZ}$ and $(\alpha_i')_{i \in \bbZ}$
if there exists $l \in \bbZ$ such that $\alpha_i = \alpha_{i +
l}'$ for all $i \in \bbZ$. The type of a characteristic sequence
$(\alpha_i)_{i \in \bbZ}$ of the above type is by definition $(0,
m)$, where $m = \min \{ l \in \bbN_+ \mid \alpha_l = \alpha_0 \}$.
In the above situation we also write $\alpha_1 \cdots \alpha_m$
instead of $(\alpha_i)_{i \in \bbZ}$.

If $(\Delta, R)$ is a gentle bound quiver, then the function
$\phi_{\Delta, R} : \bbN \times \bbN \to \bbN$, where
$\phi_{\Delta, R} (n, m)$ is the number of the characteristic
sequences of type $(n, m)$ for $(n, m) \in \bbN \times \bbN$, is a
derived invariant, i.e.\ if $(\Delta, R)$ and $(\Delta', R')$ are
derived equivalent gentle bound quivers, then $\phi_{\Delta, R} =
\phi_{\Delta', R'}$. If $\Lambda$ is the bound quiver algebra of a
gentle bound quiver $(\Delta, R)$, then we also write
$\phi_\Lambda$ instead of $\phi_{\Delta, R}$. We will write
$\phi_{\Delta, R}$ as a ``multi-set'' $[(n_1, m_1), \ldots, (n_l,
m_l)]$, where $(n, m)$ appears $\phi_{\Delta, R} (n, m)$ times.

We calculate the values of the above invariant for algebras
appearing in Theorems~\ref{theo_nondeg} and~\ref{theo_deg}, and
this will finish the proofs of these theorems. The proof of the
following lemma we leave to the reader as an easy exercise.

\begin{lemm} \label{lemm_val}
We have the following.
\begin{enumerate}

\item
If $p \in \bbN_+$ and $r \in [0, p - 1]$, then
\[
\phi_{\Lambda_0 (p, r)} = [(p, p + 2)].
\]

\item
If $p \in \bbN_+$, then
\[
\phi_{\Lambda_0' (p, 0)} = [(p + 1, p + 3)].
\]

\item
If $p_1, p_2 \in \bbN_+$, $p_3, p_4 \in \bbN$, and  $r_1 \in [0,
p_1 - 1]$, are such that $p_2 + p_3 \geq 2$ and $r_1 + p_4 \geq
1$, then
\[
\phi_{\Lambda_1 (p_1, p_2, p_3, p_4, r_1)} = [(p_1 - r_1 - 1, p_1
+ p_2), (p_2 + p_3 - 1, p_3), (r_1 + p_4, p_4)].
\]

\item
If $p_1, p_2 \in \bbN_+$, $p_3 \in \bbN$, $r_1 \in [0, p_1 - 1]$,
$r_2 \in [0, p_2 - 1]$, are such that $p_3 + r_1 + r_2 \geq 1$,
then
\[
\phi_{\Lambda_2 (p_1, p_2, p_3, r_1, r_2)} = [(p_1 - r_1 - 1,
p_1), (p_2 - r_2 - 1, p_2), (r_1 + r_2 + p_3, p_3)].
\]

\end{enumerate}
\end{lemm}

\end{document}